\documentclass[11pt, reqno]{amsart} 
\usepackage{amsthm}
\usepackage{amsmath,amssymb,latexsym,amsfonts,mathrsfs,enumerate}
\usepackage{chngcntr}
\usepackage{color}

\usepackage{mathtools}
\usepackage{fancyhdr}
\usepackage[top=30truemm, bottom=30truemm, left=25truemm, right=25truemm]{geometry}
\newcommand{\ep}{\varepsilon}
\newcommand{\nn}{\nonumber}
\newcommand{\SCR}[1]{{\mathscr #1}}

\newcommand{\CAL}[1]{{\cal #1}
}

\newcommand{\J}[1]{\left\langle #1 \right\rangle}
\newcommand{\D}[1]{{\mathscr D} \left( #1 \right)}

\newtheorem{Thm}{{\bf Theorem}}[section]

\newtheorem{Lem}[Thm]{{\bf Lemma}}
\newtheorem{Prop}[Thm]{{\bf Proposition}}
\newtheorem{Cor}[Thm]{{\bf Corollary}}
\newtheorem{Ass}[Thm]{{\bf Assumption}}

\newtheorem{Rem}[Thm]{{\bf Remark}}

\newcounter{Exami}

\def\CAL{\mathcal} 
\def\R{{\mathbb R}}
\def\N{{\mathbb N}}
\def\C{{\mathbb C}}
\def\r{\J{r}}

\def\la {\langle}
\def\ra{\rangle} 
\def\pa{\partial}
\def\ep{\epsilon}

\def\Re{{\rm Re}}
\def\Im{{\rm Im}}
\def\CAL{\mathcal} 
\def\R{{\mathbb R}}
\def\N{{\mathbb N}}
\def\C{{\mathbb C}}
\def\r{\J{r}}
\def\pa{\partial}

\def\svx{\sqrt{\varphi'(x)}}
\def\svy{\sqrt{\varphi'(y)}}
\def\svd{\sqrt{\varphi'(x) \varphi'(y)}}
\def\nn{\nonumber}

\begin{document}

\title[Embedded eigenvalues]{Absence of embedded eigenvalues for Hamiltonian with crossed magnetic and electric fields}

\author[M. Dimassi, M. Kawamoto]{Mouez Dimassi, Masaki Kawamoto}
\address{Universit\'e de Bordeaux, Institut de Math\'ematiques de Bordeaux, 351, Cours de la Lib\'eration,
33405 Talence, France}
\email{mdimassi@u-bordeaux.fr}
\address{Department of Engineering for Production, Graduate School of Science and Engineering, Ehime University, 3 Bunkyo-cho Matsuyama, Ehime, 790-8577. Japan}
\email{kawamoto.masaki.zs@ehime-u.ac.jp}
\author[V. Petkov]{Vesselin Petkov}
\address{Universit\'e de Bordeaux, Institut de Math\'ematiques de Bordeaux, 351, Cours de la Lib\'eration,
33405 Talence, France}
\email{petkov@math.u-bordeaux.fr}


\begin{abstract}

In the presence of the homogeneous electric field ${\bf E}$ and the homogeneous perpendicular magnetic field ${\bf B}$, the classical trajectory of a quantum particle on $\R^2$ moves with drift velocity $\alpha$ which is perpendicular to the electric and  magnetic fields. For such Hamiltonians the absence of the embedded eigenvalues of perturbed Hamiltonian has been conjectured.  In this paper one proves this conjecture for the perturbations $V(x, y)$ which have sufficiently small support in direction of drift velocity. 

\end{abstract}

\maketitle



\section{Introduction} 
\renewcommand{\theequation}{\arabic{section}.\arabic{equation}}
\setcounter{equation}{0}

We consider the quantum dynamics on the plane $\R^2$ in the presence of a homogeneous constant electric field
 which lies on this plane and a constant magnetic field which is  perpendiculars to this plane. Therefore the quantum system can be described by the following magnetic Stark Hamiltonian acting on $L^2(\R^2)$
\begin{align*}
H_{LS} := \frac{1}{2m} \left( 
D_X + \frac{B}{2} Y 
\right) ^2 + \frac{1}{2m} \left( 
D_Y - \frac{B}{2} X
\right) ^2  - qE \cdot {\bf X} +V, 
\end{align*}
where $D_X = -i \partial_X,\:D_Y = -i \partial_ Y$, ${\bf X} = (X,Y)  \in {\R}^2$, $m>0$, $q \neq 0$ are the position, the mass and the charge of a quantum particle and ${\bf E} = E = (E_1,E_2) \neq (0,0)$, ${\bf B} = (0,0, B)$, $B \neq 0$ stand for the electric field and the magnetic field, respectively.  Next $V$: ${\R}^2 \to \R$ is the multiplication operator by $V({\bf X})$. We assume that $V({\bf X}) $ is bounded and decays as $|{\bf X}| \to \infty$. Under some decaying conditions for the potential $V$, in \cite{DP}  and \cite{AK} it was established that 
$$\sigma (H_{LS} -V) = \sigma_{ac} (H_{LS} -V) =\R, \:\sigma_{ess} (H_{LS}) = \R.$$
Here $\sigma(L),\: \sigma_{ac}(L),\: \sigma_{ess}(L),\: \sigma_{pp}(L)$ denote the spectrum, the absolutely continuous spectrum, the essential spectrum and the point spectrum, respectively,  of the operator $L$.
In the physical literature it was conjectured that $\sigma_{pp}(H_{LS}) = \emptyset.$ This property has been proved in the following cases:
$\:\:$\\

(I) $|qE|^2 - qE \cdot \nabla V >0$ for all ${\bf X}\in{\R}^2 $ (see \cite{DP2}), \\ 

 (II) $|qE|$ is sufficiently large \cite{AK} or sufficiently small \cite{DP2}.  \\

Moreover, it was shown in \cite{DP} that\\

(III) There exists $R_0 >0$ such that $\sigma_{pp} (H_{LS}) \cap \Bigl((-\infty, -R_0 ] \cup [R_0, \infty)\Bigr) =\emptyset$  and, moreover, there exist at most a finite number of eigenvalues with finite multiplicities. \\  

In particular, (II) implies that if eigenvalues exist, then $|qE|$ is not small as well as not large. This condition seems very strange and it is natural to show that for any $|qE|$, $H_{LS}$ has no eigenvalues. The absence of point spectrum of $H_{LS}$ is an open and challenging problem. There are two major difficulties in the investigation of this problem.
If we consider the operator $H_0 = H_{LS} - V$,  first $H_0$  has double characteristics and second the electric and magnetic fields are not decreasing as $|\bf X| \to \infty$. Consequently, even for $H_0$  it is quite difficult to obtain weighted estimates for the resolvent 
$$\| \J{{\bf X}}^{-s}  (H_0 - \lambda - \pm i \ep)^{-1} \J{{\bf X}}^{- s}\|, \: s > 1/2, \: |\lambda | \gg 1$$
 uniformly with respect to $\ep \neq 0.$  In the literature there are a lot of works treating weighted resolvent estimates for the perturbations of the Laplacian. Recently the proof of suitable Carleman estimates led to several important results. We cite only some recent works \cite{V1}, \cite{V2}, \cite{V3}, \cite{GS}, where the reader may find other references. However, in these papers some decay of the potentials is assumed and this plays a crucial role in the analysis. Studying  $H_0$, we cannot  treat $H_0$ as a perturbation  of $-\Delta$ since electromagnetic
potentials do not decrease in ${\bf X}$ but increases quadratically. Usually, a  Hamiltonian with 
quadratic potential may have only bound states. Nevertheless, the presence
 of electric potential $qE \cdot  {\bf X}$ implies that  $H_0$ has only a continuous spectrum. In this direction the Hamiltonian
$H_0$  is an exceptional model in quantum physics. In the case with a  potential $V$, it is natural to consider $H_0$ as an unperturbed operator and to obtain resolvent estimates for $H_0$.   
  
  We examine the situation when the support of $V(X, Y)$ in direction of {\it drift velocity} 
  $$\alpha = (E_2 /B, -E_1/B)$$
   is sufficiently small.  Passing to new coordinates $(x, y)$, this means that the support of $V(x, y)$ has small support with respect to $y$ (see Assumption 1.1 below). We do not impose conditions on $|qE|$.\

Concerning the velocity $\alpha$, notice that according to Proposition 4.4 in Adachi-Kawamoto \cite{AK}, we have
the estimate
\begin{align*}
\mathrm{s-}\lim_{t \to \infty} 
\chi (q^2B^2 \alpha \cdot {\bf X} \leq c_1 t) \varphi (H_{LS}) e^{-itH_{LS}}   = 0, 
\end{align*} 
where $\chi$ is the characteristic function such that $\chi (s \leq a) = 1$ if $s \leq a$ and $\chi (s \leq a) = 0 $ if $s >a$, $\varphi \in C_0^{\infty} (\R)$ and  $c_1 >0$ is a suitable constant. This proposition shows that the quantum particle described by this system undergo a uniform linear motion in direction $\alpha$. By using this proposition, Kawamoto \cite{Ka} characterized the space $L^2_{pp}(H_{LS})$ of all eigenstates of $H_{LS} $, as follows: 
\begin{align*}
\psi \in L^2_{pp}(H_{LS}) \Leftrightarrow \lim_{R \to \infty} \sup_{t \in \R} \left
\|\chi (R \leq |\alpha \cdot {\bf X}| )e^{-itH_{LS}} \psi 
\right\|_{L^2(\R^2)} =0.
\end{align*}
 Hence the norms of the eigenfunctions over the region $|\alpha \cdot {\bf X}| \geq R$ goes to 0 as $R \to \infty$, it is expected that the behavior of the potential in direction perpendicular to $\alpha$ must be negligible  for the existence of eigenvalues. It is easy to see that
 \begin{align*}
 |\alpha|^{-2} (\alpha  \cdot {\bf X}) ^2 + |qE|^{-2} (qE \cdot {\bf X}) ^2 = |{\bf X}|^2, \quad \alpha \cdot qE = 0,
 \end{align*} 
 (see \S{1} of \cite{Ka}) which implies that the direction $\alpha$ is perpendicular to $qE$.\\

 In the following up to the end of the paper for simplicity we assume $m = 1/2, B = 1,q= 1.$ Introduce the change of variables 
 \begin{align} \label{eq:1.1}
 |E| x = -E \cdot {\bf X}, \qquad  |\alpha | y =  -\alpha  \cdot {\bf X} , 
 \end{align}
  hence 
 \begin{align*}
 x= - \frac{E_1 X + E_2 Y}{|E|}, \quad y = \frac{-E_2X + E_1 Y}{|E|}, \quad 
 X= - \frac{E_1 x +E_2 y}{|E|} , \quad Y = - \frac{E_2 x - E_1 y}{|E|}
 \end{align*}
 and 
 \begin{align*}
 \partial _X = - \frac{ E_1}{|E|} \partial _x - \frac{ E_2}{|E|} \partial _y , \quad 
 \partial _Y = - \frac{E_2}{|E|} \partial _x + \frac{E_1}{|E|} \partial _y .
 \end{align*}
 By using these variables, the Hamiltonian $H_{LS}$ is reduced to 
 \begin{align*}
 H_{LS} =  \left( 
 D_x + \frac{1}{2} y
 \right)^2 + \left( 
 D_y - \frac{1}{2} x
 \right)^2 + |E| x + V(x,y) 
 \end{align*}
and with the unitary transform $e^{i xy/2}$, we have 
 \begin{align*}
 e^{-i xy/2} H_{LS}e^{i xy/2} = \left( 
 D_x +  y
 \right)^2 +  D_y^2 + |E|x + V(x,y). 
 \end{align*}
The potential $V$ changes but we will denote again the new potential by $V(x,y).$
 Throughout this paper we assume $|E|=1$ and consider the reduced Hamiltonians 
\begin{align*}
H &:= H_0 + V, \\ 
H_0 &:= (D_x + y)^2 + D_y^2 + x,
\end{align*}
acting on $\D H = \D {H_0} \subset L^2(\R^2)$. In the exposition we will use the notation $\J{r} = (1 + r^2)^{1/2}, \: r = (x^2 + y^2)^{1/2}$ and similar notation for $\J{x}, 
\J{y}.$ \\
     
       The purpose of this paper is the study two problems:\\
       
      {\bf (A)} Estimates for the resolvent $\J{r}^{-\delta} (H_0 - \lambda - i\nu)^{-1}\J{r}^{-\delta}$ for $|\lambda| \gg 1,\: \nu > 0$ and $\delta > 0.$\\
      
      {\bf (B)} Absence of eigenvalues of the operator $H$.\\

The problem {\bf (A)} is examined in Section 3 and we prove the estimate 
\begin{align} \label{eq:1.2}
 \|\r^{-\delta} (H_0 -\lambda -i \nu )^{-1} \r^{-\delta}\|_{L^2 \to L^2} 
\leq C_{\delta}  \nu^{-3}  |\lambda|^{- \delta/4},
\end{align} 
where $0 < \delta \leq 2, \:0 <  \nu \leq 1$  and $|\lambda| \gg 1.$ For $\delta = 2$ we obtains the optimal decay ${\mathcal O}(|\lambda|^{-1/2})$.
The proof of Proposition 3.1 is based on  a representation of operator $e^{i t H_0}$ established in \cite{AK}. It seems that this is the first result where we have an estimate of the resolvent of $H_0$ as $|\lambda|\to \infty.$ For operators with magnetic and electric potentials having some decay a similar result with  bound  ${\mathcal O}(|\lambda|^{-1/2})$ and constant $C > 0$  uniform with respect to $\nu > 0$ has been obtain by Vodev \cite{V1}. For Stark Hamiltonian without magnetic field estimates of the resolvent are given in \cite{AD}. Concerning $H_{LS}$, it is an open problem to improve (\ref{eq:1.2})  with a constant independent of $\nu > 0$.\\

The problem {\bf (B)} is studied under the assumption

\begin{Ass}\label{Ass2}
We have $V(x, y), \partial_x V(x,y)  \in C (\R^2)$ and there exists $\eta_ 0 > 0$ such that
\begin{align} \label{eq:1.3} 
\mathrm{supp}(V) \subset \left\{ (x,y) \in \R^2 \, : \, - \eta_0 \leq y \leq \eta_0\right\}.
\end{align}
Moreover,  the potential $V(x, y) $ satisfies the estimates
\begin{align} \label{eq:1.4}
\sup_{x,y \in \R}
\J{x}^{ 2s}   |V (x,y)| \leq A_0,
\sup_{x,y \in \R}\J{x} |V_x (x,y)| \leq A_1, 
\end{align} 
with constants $ s > 1/2, \:  A_k > 0, k=0, 1$. 
\end{Ass} 
\begin{Rem}
Under Assumption $\ref{Ass2}$, it is easy to prove that the operators 
\begin{align*}
 V (H + i)^{-1}  
, \quad 
(H-i)^{-1} \left( \partial _x V \right) (H+i)^{-1} 
\end{align*} 
are compact ones. 
\end{Rem}
Our goal is to prove the absence  of embedded eigenvalues of $H$, when $V$ satisfies (\ref{eq:1.3}) with small $\eta_0.$  
To examine the non-existence of eigenvalues of $H$, first we prove in Proposition 4.2 that without {\it any assumption} on the support of $V$ there exist $R_1 > 0, R_2 > 0$ {\it independent of the support} of $V$  such that
$$\sigma_{pp}(H) \cap \Bigr( (-\infty, -R_1) \cup (R_2, \infty)\Bigr) = \emptyset.$$
This statement is more precise that the result  in \cite{DP}, where the dependence of the support of $V$ was implicit. Moreover, in contrast to \cite{DP}, the  constants $R_1, R_2$ are explicitly given and we have
\begin{equation*}
R_1 = C_1 + \|V\|_{L^{\infty}} ,\:R_2 = C_2 \|\J{x} \J{y} V_x\|_{L^{\infty}}^8, 
\end{equation*} 
where $C_1 > 0, C_2 >0$ are independent of $\eta_0$ and $V$. We see that the {\bf eigenvalues-free region} depends only of the {\it amplitudes} $A_0, \: A_1.$ 
The argument in Section 4 is based on Lemmas 2.2, 2.3 and 2.7. and we show that with a suitable weight $\varphi(x) > 0$ one has the estimates
$$\|\sqrt{\varphi(x)} (D_x + y) (H- \lambda - i)^{-1} \|_{L^2 \to L^2} \leq C \J{\lambda}^{1/2},$$
$$\|\sqrt{\varphi(x)} D_y (H- \lambda - i)^{-1} \|_{L^2 \to L^2} \leq C \J{\lambda}^{1/2},$$
with constant $C > 0$ independent of $\lambda$. In the literature such type of estimates with $\lambda = 0$ have been used without a weight $\sqrt{\varphi(x)}$. However we show in Appendix A that the operators $(D_x + y)(H-i)^{-1},\: D_y(H- i)^{-1} $ are {\it unbounded} (see Remark 2.4). We expect that the properties of these operators as well as Appendix A will be useful for further analysis.

Obviously, if $V$ satisfies Assumption 1.1 with $s \geq 3/4,$ then $V$ satisfies (\ref{eq:1.4}) with $s < 3/4$. Next in the exposition without loss of generality we assume that $1/2 < s < 3/4.$  Fixing $R = \max\{R_1, R_2\}> 0$, we establish a Mourre type
estimate for the operator $H_0$. More precisely, setting $\gamma = 2s -1 < 1/2,$ there exists a constant $C_{R, \gamma}$ independent of $\eta_0$ 
 such that
\begin{align} \label{eq:1.5}
\sup_{\lambda \in [-R,R], \nu >0} \left\| |y|^{-\gamma} F\Bigl(\frac{y}{\eta_0} \Bigr)\J{x}^{-1/2-\gamma/2} 
 (H_0-\lambda \mp i \nu )^{-1} \J{x}^{-1/2 -\gamma/2} F\Bigl(\frac{y}{\eta_0}\Bigr) |y|^{-\gamma}
\right\|_{L^2 \to L^2}  \leq C_{R, \gamma} ,
\end{align}
where  $F(t) \in C_0^{\infty}(\R: [0, 1])$ is a cut-off function such that $ F(t) = 1$ for $|t| \leq 1,$  $F(t) = 0$ for $|t| \geq 2$ (see Proposition 5.2).  It is well known that $H_0$ has no eigenvalues in $\R$, however the above estimate for the resolvent of $H_0$ is non-trivial. Since $H_0$ has only continuous spectrum, the starting point is the Mourre estimate
(\ref{eq:5.1}). Since $\gamma < 1/2,$ the weight $|y|^{-\gamma}$ is integrable around 0 and this plays an essential role.\\

Our main result is the following
\begin{Thm}\label{T1}
Let $V$ satisfy  the Assumption $\ref{Ass2}, \: \gamma = 2s - 1$ and let $C_{R, \gamma} > 0$ be the constant in $(\ref{eq:1.5}).$   Assume that
\begin{equation} \label{eq:1.6} 
\eta_0^{2\gamma} C_{R, \gamma} A_0 = c_{R, \gamma, \eta_0} < 1.
\end{equation}
Then the operator $H$ has no embedded eigenvalues. 
\end{Thm}
The condition (\ref{eq:1.6}) does not imply the smallness of the potential, it is not related to $q |E|$ as well as to the cases (I) and (II). In fact, given a potential 
satisfying Assumption \ref{Ass2}, one can choose the constant $\eta_0$ small enough in function of $A_0$ and $C_{R, \gamma}$, to obtain that there are no embedded eigenvalues of $H$.\\

Our result may be generalized to cover the case when the support of $V$ is included in a strip $\{(x, y) \in \R^2: |y - \beta| \leq \eta_0\}$ with fixed $\beta > 0$. Also, we can consider potentials having some singularities for $|x| + |y| \leq K$. These generalizations need some technical modifications, but the idea of the proof is the same. For simplicity of the exposition we are not going to treat them. 
           
Considering the case $\eta _0 \to 0$, one can choose $A_0$ large enough and such a case is closely related to the one where the potential is a delta function. When the potential $V$ is a delta function, there are interesting results due to Hauge-van Leeuwen \cite{HL}, Gyger-Martin \cite{GM} concerning the non-existence of embedded eigenvalues.   
Finally, notice that the absence of embedded eigenvalues is important for the analysis of the resonances widths (see \cite{FK1}, \cite{FK2}) .\\
 
 Setting $s = 1/2 +\gamma/2$, the proof of Theorem 1.3  is based on the equality
\begin{align*}
 |y|^{-\gamma} F\Bigl(\frac{y}{\eta_0}\Bigr)\J{x}^{-s}  (H-\lambda - i \nu)^{-1}|D_x + i|^{-1} =  |y|^{-\gamma} F\Bigl(\frac{y}{\eta_0}\Bigr) \J{x}^{-s}(H_0 -\lambda - i\nu)^{-1}|D_x + i|^{-1} \\
 -  |y|^{-\gamma} F\Bigl(\frac{y}{\eta_0}\Bigr)\J{x}^{-s}  (H_0 -\lambda - i \nu )^{-1}  V (H-\lambda - i \nu)^{-1} |D_x + i|^{-1} , \nu > 0,\: 0 < \gamma < 1/2,\end{align*}
 where $\lambda \in \R$ and $F\Bigl(\frac{y}{\eta_0}\Bigr)$ is a cut-off function equal to 1 for $| y| \leq \eta_0.$
 
 By the condition of $V$, we have $V =  |y|^{-2\gamma} F({y}/{\eta_0}) \cdot |y|^{2\gamma} V$ and for $\eta_0^{2\gamma}C_{R, \gamma} \|\J{x}^{1 + \gamma} V\|_{L^{\infty}}= c_{R,\gamma, \eta_0} < 1$  we may estimate
$$\sup_{\lambda \in [-R, R], \nu > 0} \Bigl\| |y|^{-\gamma} F\Bigl(\frac{y}{\eta_0}\Bigr) \J{x}^{-s}(H-\lambda - i\nu)^{-1} |D_x + i|^{-1}  \Bigr\|_{L^2 \to L^2}  \leq \frac{B_{R, \gamma}}{1 - c_{R, \gamma, \eta_0}}.$$
If $\psi \in L^2(\R^2)$ is an eigenfunction of $H$ with eigenvalues $\lambda \in [-R, R]$, we show in the Appendix B that $D_x \psi \in L^2$ and one obtains easily a contradiction with the above estimate.
 Following  this approach, one needs to establish uniform estimate (\ref{eq:1.5}). To cover more general cases of potentials, it is necessary to obtain estimate similar to (\ref{eq:1.5})  with more general weights and this is an interesting open problem.   
 
 The plan of the paper is as follows. In Section 2 we prove some preliminary results including Lemma 2.3, Lemma 2.7 and Proposition 2.9 which are used in the next sections. In Section 3 we examine the estimates of the resolvent of $H_0$.  The absence of large eigenvalues of $H$ is studied in Section 4. Mourre type estimates are proved in Section 5 and Theorem 1.3 is established in Section 6. In the Appendix A and B we prove some technical results. Finally, notice that we use the Assumption 1.1 for the support of $V(x, y)$ with respect to $y$ only in Section 6. The results in other sections concerning $H$ hold without any restriction on the support of $V$.
 
\section{Preliminaries} 
\renewcommand{\theequation}{\arabic{section}.\arabic{equation}}
\setcounter{equation}{0}

In this section we prove some lemmas which are necessary for the exposition. Throughout this section to the end of this paper, $ \left\|  \cdot \right\|_{L^2(\R^2)}$ and $ \left\|  \cdot \right\|_{\SCR{B}(L^2(\R^2))}$ are denoted as $\| \cdot \|$ and $(\cdot, \cdot)$ denotes the inner product on $L^2(\R^2)$. Also we denote by $\D{A}$ the domain of the operator $A$. We write $r := \sqrt{x^2 +y^2}$ and $\J{\cdot} = (1 + \cdot ^2)^{1/2}$.

\begin{Lem}[Interpolation Theorem]\label{LL1}
Let $A$ and $B$ be positive selfadjoint operators on $L^2(\R^2)$ and let $T$ be a bounded operator on $L^2(\R^2).$  Assume that with constants $\alpha _0, \beta_0, \alpha_1 , \beta _1 \geq 0$ and $C_0, C_1 >0$ we have
\begin{align*}
& \left\| 
A^{\alpha _0} T B^{\beta _0}
\right\| \leq C_0, \\ 
& \left\| 
A^{\alpha _1} T B^{\beta _1}
\right\| \leq C_1.
\end{align*}  
Then for all $0<\theta' <1$, setting $\alpha _{\theta'} = \alpha _0 (1- \theta') + \alpha _1 \theta'$ and $\beta _{\theta'} = \beta _0 (1- \theta') + \beta _1 \theta'$, one has
\begin{align*}
\left\| 
A^{\alpha _{\theta'}} T B^{\beta _{\theta'}}
\right\| \leq C_0^{1- \theta'} C_1^{\theta'}.
\end{align*}

\end{Lem}
\begin{proof}
We can find the proof of this lemma for example  in \S{6} in Isozaki \cite{Is}. We will give a sketch of the proof based on the Hadamard's three line theorem. Recall this theorem.\\ 
Let $f(z)$ be an analytic function on $\Omega_0 := \left\{ z =
x+iy:\:  0<x<1 , -\infty < y < \infty
\right\}$ which is bounded on $\overline{\Omega_0} := \left\{ 
z =x+iy :\,  \, 0 \leq x \leq 1 , -\infty < y < \infty
\right\}$. Then if one has
\begin{align*}
& \sup_{- \infty < y < \infty} |f(iy)| \leq M_0, \\ 
& \sup_{- \infty < y < \infty} |f(1+ iy)| \leq M_1,
\end{align*}
 then for all $0<x<1$ we have the estimate 
\begin{align*}
\sup_{- \infty < y < \infty} |f(x + iy)| \leq M_0^{1-x} M_1 ^x. 
\end{align*}

 Let $u,v \in L^2(\R^2)$, and let $E_A(\cdot),\: E_B(\cdot)$ be the  spectral decompositions of $A$ and $B$, respectively. Let $I\subset \R$ be some bounded interval. Define 
\begin{align*}
f(z) := \left( 
E_A (I) A^{\alpha _z} TB^{\beta _z} E_B(I) u,v 
\right),   
\end{align*}
where $\alpha _z:= (\alpha _1 - \alpha _0) z + \alpha _0$ and $\beta_z: = (\beta _1- \beta _0)z + \beta _0$. Then $f(z)$ is an analytic function on $0< \mathrm{Re}z <1  $ which is continuous and bounded on $0 \leq \mathrm{Re} z \leq 1$ and for all $y \in \R$ one gives
\begin{align*}
|f(iy)| \leq C_0  \| u \|\| v \|, \quad |f(1+ iy)| \leq C_1 \| u \| \| v \|.
\end{align*}
 By the theorem above, for all $0< \theta' <1$, we get 
\begin{align*}
|f(\theta')| \leq C_0^{1- \theta'} C_1^{\theta'} \| u\| \|v \|. 
\end{align*}
Since $C_0,C_1$ are independent of $I$, by taking $I \to \R$, one completes the  proof of Lemma \ref{LL1}.   \end{proof}

Introduce a positive function $\rho(x) \in C^{\infty} (\R)$ such that for some fixed $a > 1$ we have $\rho(x) = - \frac{1}{x},\:x \leq -a$, $\rho(x) = 2x, 2x > a.$
 \begin{Lem} For $f \in \D{H}$ we have
 \begin{equation} \label{eq:2.1}
 \iint_{\R ^2} \rho(x) \Bigl( |(D_x + y)f|^2 + |D_y f|^2\Bigr) dxdy \leq C (\|H f\|^2 +\|f\|^2).
 \end{equation}
\end{Lem} 
\begin{proof} Consider for $f \in C_0^{\infty} (\R^2) $ the integral
$$0 \leq \iint_{\R ^2} |(H_0 - \rho) f |^2 dxdy = \iint_{\R^2} \Bigr( |H_0 f|^2 + (\rho^2 -2 x \rho + \rho'') |f|^2\Bigr) dxdy$$
$$ - 2\iint_{\R ^2} \rho \Bigl( |(D_x +  y)f|^2 + |D_y f|^2\Bigr) dxdy.$$
Here we have used that by integration by parts one obtains
$$ -\mathrm{ Re}\:\iint_{\R ^2}\Bigl[ ((D_x +  y)^2 f )\rho \bar{f} +\rho f \overline{(D_x + y)^2 f}\Bigr] dxdy = \iint_{\R^2} \Bigl( - 2\rho |D_x + y) f|^2 + \rho'' |f|^2\Bigr) dxdy$$
and similarly one transforms the integral with $D_y^2 f.$ Clearly with a constant $C_0 > 0 ,$ one has
$$\rho^2 - 2 x \rho + \rho'' < C_0,$$
hence
$$ 2\iint_{\R ^2} \rho(x) \Bigl( |(D_x +  y)f|^2 + |D_y f|^2\Bigr) dxdy \leq C_1( \|H_0 f\|^2 +  \|f\|^2) \leq C_2( \|H f\|^2 +\|f\|^2).$$
Since $H$ is a closed operator, for every $f \in \D{H}$ there exists a sequence of functions $f_n \in C_0^{\infty} (\R^2)$ such that $f_n \to f,\: Hf_n \to H f$ in $L^2.$ Taking the limit $n \to \infty$, we obtain the result. \end{proof}
The above Lemma is analogous to Lemma 1 in \cite{Va} for Stark Hamiltonian.\\

Consider a function $0 < \varphi(x) \leq A$ defined by
$$\varphi(x) = \begin{cases} 
 \tan^{-1} (x-a + \pi/4)  ,& x \geq a, \\ 
\varphi_1(x), &  -a < x < a\\
-\frac{1}{ x} , & x \leq -a, 
\end{cases}$$
where $\frac{1}{a} \leq \varphi_1 (x) \leq 1,\: |x|\leq a$ is a smooth function so that $\varphi(x ) \in C^3 ({\R})$ and $\varphi'(x) >0$ for all $x \in \R$.
 \begin{Lem} We have the estimates
 \begin{equation} \label{eq:2.2}
 \Bigl\| \varphi^{1/2} (x)(D_x + y) (H - \lambda - i)^{-1} \Bigr\| _{L^2 \to L^2} \leq C \J{\lambda}^{1/2},
 \end{equation} 
  \begin{equation} \label{eq:2.3}
 \Bigl\| \varphi^{1/2}(x) D_y (H - \lambda - i)^{-1} \Bigr\| _{L^2 \to L^2} \leq C \J{\lambda}^{1/2}
 \end{equation}  
 with a constant $C = C_a > 0$ independent of $\lambda.$
\end{Lem}
\begin{proof}  By using the resolvent equality
$$(H - \lambda - i)^{-1} = (H_0 -\lambda - i)^{-1} - (H_0 - \lambda - i) ^{-1} V( H- \lambda - i)^{-1},$$
it is sufficient to prove the estimates with $H$ replaced by $H_0.$ We apply the unitary operator $e^{i \lambda D_x}$ giving a shift $x \rightarrow x + \lambda,$ and obtain
$$e^{i \lambda D_x} \varphi^{1/2}(x)  (D_x + y)e^{- i \lambda D_x} e^{i \lambda D_x} (H_0 - \lambda - i)^{-1} e^{- i \lambda D_x}$$
$$= \varphi^{1/2} (x + \lambda) \rho^{-1/2}(x)  \Bigl(\rho^{1/2}(x) (D_x + y)(H_0 - i)^{-1} \Bigr).$$
On the other hand,
$$\Bigl \|\rho^{1/2}(x) (D_x + y)(H_0 - i)^{-1} \Bigr\|_{L^2 \to L^2}  \leq C.$$
In fact, we apply (\ref{eq:2.1}) replacing $H$  by $H_0$ and choose $ f = (H_0 - i)^{-1} g$. This yields
$$\|\rho^{1/2} (x)(D_x + y) (H_0 - i)^{-1} g\|^2 \leq C\Bigl( \|H_0 (H_0 - i)^{-1} g\|^2 + \|(H_0 - i)^{-1} g\|^2\Bigr)\leq C_1 \|g\|^2.$$
It remains to prove the estimate 
\begin{equation} \label{eq:2.4}
\varphi(x + \lambda) \rho^{-1}(x) \leq C_2 ( 1 + |\lambda|)
\end{equation}
with $C_2= C_2(a)  > 0$ independent of $\lambda.$ 
For $ x \geq -2 a$ the function $\rho^{-1} (x)$ is bounded by a constant $B_a$ depending on $a$ and $\varphi(x + \lambda) \leq A$, hence we have (\ref{eq:2.4}) . We are going to study the case $x < -2 a.$ We have three subcases: (i) $ |\lambda| \leq a,\:  x < -2 a$, (ii) $|\lambda|  > a, \: - 2 |\lambda| \leq x < -2a$, (iii) $|\lambda| > a, \: x < - 2 |\lambda|.$
Clearly, in the subcase (i) one has  $x + \lambda\leq -2 a +|\lambda| \leq -a$ and $\varphi(x + \lambda) \rho^{-1} (x) = \frac{2 |x|}{|x+ \lambda|} \leq 2 + \frac{2|\lambda|}{a}.$
In the subcase (ii) we have $\rho^{-1} (x) = |x| \leq 2|\lambda|.$ In the subcase (iii) we have $x + \lambda < - |\lambda| \leq - a$ and $\varphi(x+ \lambda)\rho^{-1}(x)   = \frac{|x|}{|x+\lambda|} \leq 1 +\frac{|\lambda|} {a}.$ Thus we obtain (\ref{eq:2.4}).
For $D_y$ we apply the same argument.
\end{proof}
\begin{Rem} The presence of the factor $\varphi^{1/2}(x)$ in the estimates $(\ref{eq:2.2})$, $(\ref{eq:2.3})$ is important for the boundedness of these operators. In fact, the domain
$\D{H}$ is not included in the domains $\D{ D_x + y}, \: \D{D_y}$ and both operators $(D_x + y)(H- i)^{-1},\: D_y(H- i)^{-1}$ are unbounded. We prove this property in Appendix A.
\end{Rem}
\begin{Rem} It is clear that the estimates $(\ref{eq:2.2})$, $(\ref{eq:2.3})$ hold with $\varphi(x)$ replaced by $\J{x}^{-1} .$
\end{Rem}

\begin{Cor}
For every $0 < \gamma \leq 1$ we have the estimate
 \begin{equation} \label{eq:2.4bis}
 \Bigl\| \J{x}^{-\gamma/2} \J{D_y}^{\gamma}  (H - \lambda - i)^{-1} \Bigr\| _{L^2 \to L^2} \leq C_{\gamma}  \J{\lambda}^{\gamma/2}
 \end{equation}  
\end{Cor}
\begin{proof} Writing $\J{D_y}  = \frac{ \J{D_y}}{D_y + i} (D_y + i),$ we deduce that (\ref{eq:2.4bis}) holds with $\gamma = 1.$ Next we apply the interpolation Lemma 2.1 between
$$\|\J{x}^0  \J{D_y}^0 \J{x}^0(H - \lambda - i)^{-1} \| \leq 1$$
 and 
$$\Bigl\|\J{x}^{-1/4} \J{D_y}\J{x}^{-1/4}  (H- \lambda - i)^{-1}\Bigr\| \leq C\J{\lambda}^{1/2}.$$

\end{proof} 
Notice that we have
$$\varphi'(x) = \begin{cases}  \Bigl( 1+  (x - a + \pi/4)^{2}\Bigr)^{-1}  ,\: x \geq a,\\
\varphi_1'(x),\: - a \leq x < a,\\
x^{-2} , \: x < -a,\end{cases}$$
which implies
\begin{align*}
(\varphi'(x))^{1/2}  \leq C_2 \sqrt{\varphi(x)},
\end{align*}
hence the estimates (\ref{eq:2.2}), (\ref{eq:2.3}) hold with $\varphi(x)$ replaced by $\varphi'(x).$

For the eigenfunctions of $H$ we need a more precise result.

\begin{Lem}
Let $\psi$ and $\lambda$ be an eigenfunction and eigenvalue of $H$. Moreover, suppose that $\| \psi \| = 1$.
Then  we have the estimates
\begin{align}\label{eq:2.5}
\left\| 
\sqrt{\varphi'(x)} (D_x +y) \psi
\right\| \leq C \J{\lambda}  ^{1/4},
\end{align}

\begin{align}\label{eq:2.6}
\left\| 
\sqrt{\varphi(x)} \sqrt{\varphi'(y)}  D_y \psi
\right\| \leq C \J{\lambda} ^{3/8}
\end{align}
with $C = C_a > 0$ independent of $\psi$ and $\lambda$ and the support of $V$.
\end{Lem}
\begin{proof}
   By a direct calculus we obtain a representation for the  commutator 
\begin{align*}
0 &= \Bigl( 
i \left[ H, \varphi(x) (D_x + y) + (D_x + y) \varphi(x)  \right] \psi , \psi
\Bigr) \\ &= \Bigl( \left( 
4(D_x + y) \varphi'(x) (D_x +y) + 4 \varphi(x)  D_y -\varphi(x)( 1 + V_x) + 2 \varphi'''(x)
\right)  \psi , \psi\Bigr). 
\end{align*}
Hence we have 
\begin{align*}
\left\| 
\sqrt{\varphi'(x) } (D_x +y) \psi
\right\|^2 \leq C_1  + \left\| \varphi(x) D_y \psi \right\|,
\end{align*}
where the constant $C_1 > 0$ depends only on $\varphi(x)$ and $\|V_x\|_{L^{\infty}}.$ 

Applying Lemma 2.3, one deduces
$$\|\varphi(x) D_y \psi\| \leq C \sqrt{A} \|\varphi^{1/2}(x)  D_y (H- \lambda - i)^{-1}\psi\| \leq C_2 \J{\lambda}^{1/2} $$
and 
\begin{align*}
\left\| 
\sqrt{\varphi'(x)} (D_x +y) \psi
\right\|^2 & \leq C_3 \J{\lambda}^{1/2}.
\end{align*}

Now we pass to the analysis of the estimate containing $D_y$.  In a  similar way one has
\begin{align*}
0 &= \Bigl( i \left[ H,  D_y\varphi(x)\varphi(y) + \varphi(y)\varphi(x) D_y  \right] \psi , \psi
\Bigr) \\ &=  \Bigl( ( 4 D_y \varphi(x)\varphi'(y) D_y   -2\varphi(x)\varphi'''(y)) \psi, \psi\Bigr)\\&
\quad +4 \mathrm{Re}\Bigl ( (D_x +y) \varphi'(x) \varphi(y) D_y\psi, \psi\Bigr)  -2 \mathrm{Im} \left( 
D_y \varphi (y) \varphi (x) \psi , V \psi
\right) \\& \qquad  -2 \mathrm{Im} \left( 
\varphi (x) \varphi (y) D_y \psi , V \psi
\right) \\ & \geq 
4 \left\| \sqrt{\varphi(x) \varphi'(y)} D_y \psi \right\|^2 \\ & \quad - 4 A\left\| 
\sqrt{\varphi'(x) } (D_x +y) \psi\right\| \left\|\sqrt{\varphi'(x)} D_y \psi\right \|\\ & \qquad
- C\Bigl( 1 +  \left\|\sqrt{\varphi(x)} D_y \psi\right\|\Bigr).
\end{align*}
Applying the estimates
$$
\left\| \sqrt{\varphi'(x) } (D_x +y) \psi
\right\| \leq C_a \J{\lambda} ^{1/4},$$
$$ \quad \left\| \sqrt{\varphi'(x) } D_y \psi\right\| \leq C_a\|\J{x}^{-1} D_y \psi\| \leq C_a \|\sqrt{\varphi (x)} D_y \psi\| \leq C_a  \J{\lambda} ^{1/2},$$
we obtain the result. \end{proof}

It is obvious that the estimates (\ref{eq:2.5}) and (\ref{eq:2.6}) hold with $\sqrt{ \varphi'(x)} $ and $\sqrt{\varphi'(y)}$ replaced by $\J{x}^{-1}$ and $\J{y}^{-1}$, respectively.

\begin{Rem}
Notice that by Lemma $2.2$, we obtain 
$$\J{r}^{-1/2} (D_x + y) \varphi \in L^2, \:\J{r}^{-1/2} D_y \varphi \in L^2,\: \varphi \in \D{H}.$$
Let $H_{ev}$ be the space generated by the eigenfunctions of $H$.  By the  closed graph theorem the operators $\J{{r}}^{-1/2} (D_x + y)$ and $\J{{r}}^{-1/2} D_y$ are bounded as operators from $H_{ev}$ to $L^2( \R^2).$ 
Therefore for every eigenfunction $\psi$ of $H$ we have the estimate
\begin{equation} \label{eq:3.10}
\|\J{{r}}^{-1/2} (D_x + y) \psi\| + \|\J{{r}}^{-1/2} D_y \psi \| \leq B \|\psi\|,
\end{equation} 
where $B > 0$ is  independent of $\psi$. However,  the constant $B$ in general could depend on the support of $V$.
\end{Rem}

Let $F(t) \in C_0^{\infty}(\R)$ be a function such that $0 \leq F(t) \leq 1,\: F(t) = 1$ for $|t| \leq 1,$  $F(t) = 0$ for $|t| \geq 2.$ For $ c > 0$ define $F_c(t) = F(\frac{t}{c}).$\begin{Prop} \label{P1}
Let $0 < \gamma < 1/2, \: \beta \in \R $ and $0 < \eta _0 < 1$. Then the operator
\begin{align*}
|y- \beta|^{- \gamma } F_{ \eta _0} (y- \beta) {\J{x}^{-\gamma/2}} \J{H}^{-\gamma}
\end{align*} 
is bounded and its bound is independent of $\eta_0$ and $\beta$.
\end{Prop}
\begin{proof}
First consider the case  $\beta = 0$. Let $f(t) \in C_0^{\infty}(\R)$ be a function such that $f(t) = 1$ for $|t| \leq 2.$ Then $f(t) F_{\eta_0}(t) = F_{\eta_0}(t)$ and for any $\phi \in C_0^{\infty} (\R^2)$, it is enough to prove that 
\begin{align*}
\CAL{A} := \left\| 
|y|^{- \gamma } f(y) {\J{x}^{-\gamma/2}}  \J{H}^{-\gamma} \phi 
\right\| \leq C_{\gamma} \| \phi \|
\end{align*} 
with a constant $C_{\gamma}$ dependent of $\gamma$ but independent of $\eta_0.$ By  simple calculation  one has
\begin{align*}
\CAL{A} &\leq  \left\| 
|y|^{- \gamma } f(y)  {\J{x}^{-\gamma/2}} \J{H}^{-\gamma} \phi 
\right\|\\
 &  \leq \left\| |y|^{-\gamma} \J{D_y}^{- \gamma}  \right\|_{L^2(\R_y)\to L^2(\R_y)} \left\| \J{D_y}^{\gamma} f(y) {\J{x}^{-\gamma/2}}  \J{H}^{-\gamma} \phi  \right\|. 
\end{align*}
Here we apply the fractional Sobolev inequality (see, e.g., Stein-Weiss \cite{SW} {and Yafaev \cite{Ya}}) 
\begin{align*}
\left\| 
|y|^{-\gamma} u
\right\|_{L^2({\R}_y)} \leq C_{\gamma} \left\| (D_y^2)^{\gamma/2} u  \right\|_{L^2({\R_y})} 
\end{align*}
for $u \in \D{(D_y^2)^{\gamma/2}} \subset L^2({\R^2})$. Then $\CAL{A}$ can be estimated by
\begin{align*}
\CAL{A} & \leq 
  C_{\gamma}  \|  \J{D_y} ^{\gamma} f(y) {\J{x}^{-\gamma/2}} \J{H}^{-\gamma} \phi \| .
\end{align*}
On the other hand, the norms of the operators 
\begin{align*}
&\left\| 
\J{D_y} ^0 \J{x}^0 \cdot f(y)\cdot  \J{H}^{-0}
\right\|_{L^2 \to L^2} , \\ & \left\| 
\J{D_y}^1 \J{x}^{-1/2} \cdot f(y) \cdot  \J{H}^{-1}
\right\|_{L^2 \to L^2}  \leq \left\| {\J{x}^{-1/2}}\J{D_y} f(y) (H-i)^{-1} \cdot (H -i) \J{H}^{-1}
\right\|_{L^2 \to L^2}
\end{align*}
are bounded. In fact, we write $(H- i)^{-1}  = (H_0 - i)^{-1} - (H_0 - i)^{-1}V (H- i)^{-1}$ and one applies Corollary 2.6 to estimate $\la x \ra^{-1/2} D_y(H_0- i)^{-1} $. Since $\J{x}^{-1/4} \J{D_y} \J{x}^{-1/4}$ is selfadjoint, by using the interpolation Lemma \ref{LL1}, we conclude that
\begin{align*} 
\left\|\J{D_y}^{\gamma}\J{x}^{-\gamma/2}   \cdot f(y) \cdot  \J{H}^{-\gamma}
\right\|_{L^2 \to L^2} 
\end{align*}
is bounded. Therefore
$\CAL{A} \leq  C_{\gamma} \| \phi \|$ 
 and we obtain the estimate. \\

 Now consider the case when $\beta \neq 0$. We have
  \begin{align*}
& \left\| |y-{\beta}|^{- \gamma} f(y- \beta) {\J{x}^{-\gamma/2}} \J{H}^{-\gamma} \right\| _{L^2 \to L^2}  \\ & \leq 
\left\||y-{\beta}|^{- \gamma} \J{D_y}^{- \gamma} \right\|_{\SCR{B}(L^2(\R _y))}  \left\|  \J{D_y}^{\gamma} f (y- \beta) {\J{x}^{-\gamma/2}} \J{H} ^{-\gamma}\right\|_{L^2 \to L^2} 
 \\ & \leq 
\left\|e^{-i{\beta}D_y} |y|^{- \gamma} \J{D_y}^{- \gamma} e^{i{\beta}D_y} \right\|_{\SCR{B}(L^2(\R _y))}   \left\| \J{D_y}^{\gamma} f (y- \beta) {\J{x}^{-\gamma/2}} \J{H} ^{-\gamma}  \right\|_{L^2 \to L^2}  \leq C, 
\end{align*}
noting that $\sup_ y |f'(y - \beta)|$ is independent of $\beta.$  \end{proof}

\section{Estimates of the resolvent of $H_0$}
\renewcommand{\theequation}{\arabic{section}.\arabic{equation}}
\setcounter{equation}{0}

 In the section we establish a decay estimate for
$$\|f (H_0 - \lambda -i \nu)^{-1} g\|_{L^2 \to L^2} $$
with $\nu > 0$ and $|\lambda|\to \infty.$ In \cite{FK2} the case when $f, g \in C_0^{\infty}(\R^2)$ has been studied,
 while in \cite{AK} the situation with $f, g \in L^p(\R^2),\: p >2$ was examined. We prove the following more precise result which has an independent interest.
\begin{Prop}\label{PP1}
 Consider the operator
\begin{align*}
M_{\delta}(\lambda, \nu) := \r^{-\delta} (H_0 -\lambda -i \nu )^{-1} \r^{-\delta}, \: 0 < \delta \leq 2,
\end{align*}
where $\lambda \in \R,\: 0 < \nu \leq 1$. Then for $0 < \theta \leq 1/2$ and $ |\lambda| \geq 1$, there exists a constant $C = C(\theta)>0$ such that 
\begin{align} \label{aX0}
 \left\|M_{\delta}(\lambda, \nu) \right\|_{L^2 \to L^2}  \leq C \nu^{-1}  \left( |\lambda|^{- \theta} + (1 +\nu) |\lambda|^{ -1} +( 1 + \nu^{-2})  |\lambda|^{ \theta -1}
\right)^{\delta /2}.
\end{align}

\end{Prop}
\begin{proof}
We consider only the case $\lambda >0$, since for $\lambda < 0$ the proof is similar. Set 
\begin{align*}
\omega_n := \left\{ 
t \in [0,\infty) \, : \, |t -n \pi | \leq \lambda ^{- \theta}
\right\}, \quad n \in \N \cup \{ 0 \} , 
\end{align*} 
and $\Omega = \bigcup_{n=0}^{\infty} \omega_n$.  
By the integral formula for the resolvent, we have 
\begin{align*} 
M_{2}(\lambda, \nu) &= i \J{{r}}^{-2} \int_0^{\infty} e^{-it(H_0 - \lambda -i \nu)} \J{{r}}^{-2} dt  = 
K_1 + K_2 
\end{align*}
with 
\begin{align*}
K_1 := i \J{{r}{}}^{-2} \int_{\Omega }  e^{-it(H_0 - \lambda -i \nu)} \J{{r}{}}^{-2} dt  
\end{align*} 
and 
\begin{align*}
K_2 &:= i\J{{r}{}}^{-2} \int_{[0, \infty) \setminus \Omega  }  e^{-it(H_0 - \lambda -i \nu)} \J{{r}{}}^{-2} dt
\\ &= i \J{r}^{-2}  \int_{ [0, \infty) \setminus \Omega  }  \frac{1}{i \lambda} \left( 
\frac{d}{dt} e^{it \lambda}
\right) e^{-it(H_0  -i \nu)} \J{{r}{}}^{-2} dt
  \\ &= 
\frac{1}{\lambda} \sum_{n=0}^{\infty}\left[ 
\J{{r}}^{-2} e^{-it(H_0 - \lambda -i \nu)} \J{{r}}^{-2}
\right]\big |_{t =n\pi + \lambda^{-\theta}}^{t = (n+1) \pi - \lambda^{-\theta}} \\ & \quad + \frac{i}{\lambda} \int_{[0, \infty)  \setminus \Omega} \J{{r}}^{-2} e^{-it (H_0 - \lambda -i \nu)} (H_0 -i \nu) \J{{r}}^{-2} dt.
\end{align*} 
For $\phi \in L^2$ one gets
\begin{align} \nn
\left\| 
K_1 \phi \right\| & \leq \sum_{n=0}^{\infty} \left\| 
\int_{\omega_n}  \J{{r}}^{-2} e^{-it (H_0 - \lambda -i \nu)} \J{{r}}^{-2} \phi dt
\right\|_{L^2(\R ^2)} \\ &  \nn  \leq 
\sum_{n=0}^{\infty} 
\int_{\omega_n} \left\|  \J{{r}}^{-2} e^{-it (H_0 - \lambda -i \nu)} \J{{r}}^{-2} \phi  
\right\|_{L^2(\R^2)} dt  
 \\ & \nn  
\leq
\sum_{n=0}^{\infty} \int_{n\pi-\lambda ^{- \theta}}^{n \pi + \lambda ^{- \theta}} \left\|   \J{{r}}^{-2} \right\|^2_{L^{\infty}}  \| \phi \|_{L^2(\R^2)}e^{-t \nu} dt \\ &= 
\frac{  \|\r^{-2}\|^2_{L^{\infty}} \| \phi \|_{L^2(\R^2)} } {\nu} \left( 
e^{\nu \lambda ^{- \theta}} - e^{-\nu \lambda ^{- \theta} }
\right) \sum_{n=0}^{\infty} e^{-n \nu \pi} \leq 2C\frac{ e^{\nu \lambda^{-\theta}}}{1 - e^{-\nu \pi}} \lambda^{- \theta}   \| \phi \|_{L^2(\R^2)}   \nn \\ & \leq C_1  \nu^{-1} \lambda^{-\theta}    \| \phi \|_{L^2(\R^2)}  \label{aX1}.
\end{align}
Here for $0 < \nu \leq 1$ we have used the elementary inequality $1 - e^{-\nu \pi}\geq  \frac{ \pi \nu}{2e ^{\pi}}.$ Next
\begin{align*} 
& \frac{1}{\lambda} \left\|\sum_{n=0}^{\infty}\left[ 
\J{{r}}^{-2} e^{-it(H_0 - \lambda -i \nu)} \J{{r}}^{-2}
\right]\big |_{t =n\pi + \lambda^{-\theta}}^{t = (n+1) \pi - \lambda^{-\theta}} \phi\right\| \\ 
& \leq  \lambda^{-1}  \|\r^{-2}\|^2_{L^{\infty}} \| \phi \|_{L^2(\R^2)}  \sum_{n=0}^{\infty} \left( 
e^{-\nu ( (n+1) \pi - \lambda^{- \theta}) } + e^{-\nu (n \pi + \lambda^{- \theta})}
\right) \\ & \leq 
2 C \lambda^{-1} e^{\nu \lambda^{- \theta}} \|\phi\|_{L^2} \sum_{n=0}^{\infty} e^{- \nu n \pi} \leq C\frac{2 e^{\nu \lambda^{-\theta}}}{1 - e^{-\nu \pi}} \lambda^{-1}\| \|\phi\|_{L^2}   \leq C \nu^{-1} \lambda^{-1}\| \|\phi\|_{L^2}
\end{align*}
and
\begin{align*}
\frac{1}{\lambda} \left\| 
\int_{[0, \infty)  \setminus \Omega} \J{{r}}^{-2} e^{-it (H_0 - \lambda -i \nu)} \nu \J{{r}}^{-2} dt
\right\|_{L^2}  \leq \frac{\nu}{\lambda} \|\r^{-2}\|^2 \int_0^{\infty} e^{- \nu t} dt =  \lambda^{-1}\|\r^{-2}\| ^2.
\end{align*}
Therefore,
\begin{align} \label{aX2}
\left\| K_2 \right\|_{L^2 \to L^2}  \leq C \nu^{-1} \lambda^{- 1}  \| \r^{-2} \|^2  
+ \lambda ^{-1} \|K_3 \|_{L^2 \to L^2} 
\end{align}
with 
\begin{align*}
K_3 =  \int_{[0,\infty) \setminus \Omega } \J{{r}}^{-2} H_0 e^{-it(H_0 - \lambda -i \nu)} \J{{r}}^{-2} dt. 
\end{align*}

Now we estimate $\|K_3\|_{L^2 \to L^2} $. Let ${\bf x} = (x,y)$. By an application of the formula (4.6) in  \cite{AK}, we have the following representation of the operator $e^{-i t H_0}$ (For the operator $H_{LS}$ in \cite{AK} one chooses the constants $q=B=1$, $m=1/2$, $\omega = 2$, $E_1 =-1,\: E_2 =0$, $\nu=0$, $\tilde{\nu} = \omega = 2$, $\theta = \pi$ and $E_0 =1$)
\begin{align*}
& \left( e^{-itH_0} \phi  \right)({\bf x})= (e^{-i xy/2} e^{-itH_{LS}} e^{i xy/2} \phi)({\bf x}) \\
&= 
\frac{1}{4 \pi i \sin (t)} \int_{{\R}^2} e^{-ia(t)}e^{-ix y/2} e^{ib(t) \cdot {\bf x}} e^{-i c(t) \cdot {\bf A} ({\bf x})} e^{-i {\bf w} \cdot {\bf A} ({\bf x}-c(t))  } e^{i (\cot t) ({\bf x}-c(t) -{\bf w})^2/4 } e^{iw_1w_2 /2} \phi ({\bf w}) d {\bf w} \\ 
&=: \frac{1}{4 \pi i \sin (t)} \int_{{\R}^2} K(t,{\bf x}, {\bf w}) \phi ({\bf w}) d {\bf w} 
\end{align*}
with ${\bf w} = (w_1,w_2) \in {\R}^2$,
\begin{align*}
{\bf A} ({\bf x}) = (-y/2,x/2), \quad a(t) = \int_0^t \left( 
{b(s)^2} + 2 b(s) \cdot {\bf A} (c(s))
\right) ds
\end{align*}
and $b(t) = (b_1(t),b_2(t)),\: c(t) = (c_1(t),c_2 (t))$ with 
\begin{align*}
b_1(t) = - (\sin (2t) )/2, \quad b_2 (t) = (1-\cos (2t)) /2,
\end{align*}
\begin{align*}
c_1 (t) = \cos (2t), \quad c_2 (t) = t - \sin (2t). 
\end{align*}
Simple calculation shows that 
\begin{align*}
& y\partial _x K(t,{\bf x}, {\bf w}) = y\left( - i \frac{y}{2} + 
ib_1(t) - i \frac{c_2(t)}{ 2} -i \frac{w_2}{2} +i\frac{\cot t}{2} \left( 
x -c_1 (t) - w_1
\right)\right) K(t, {\bf x} , {\bf w})
 , \\ 
& \partial _x^2K(t,{\bf x}, {\bf w}) = \left( \left( - i \frac{y}{2} +
ib_1(t) - i \frac{c_2(t)}{ 2} -i \frac{w_2}{2} +i \frac{\cot t}{2} \left( 
x -c_1 (t) - w_1 \right)
\right) ^2  + i\frac{\cot t}{2}  \right)  K(t, {\bf x} , {\bf w})
\end{align*}
and 
\begin{align*}
& \partial _y^2 K(t,{\bf x}, {\bf w})  \\ 
&= 
 \left( \left(- i \frac{x}{2}   + 
ib_2(t) + i \frac{c_1(t)}{ 2} + i \frac{w_1}{2} + i\frac{\cot t}{2} \left( 
y -c_2 (t) - w_2 \right)
\right) ^2  + i\frac{\cot t}{2}  \right)  K(t, {\bf x} , {\bf w}).
\end{align*} 
Thus we deduce 
\begin{align*}
& \J{r}^{-2} H_0 K(t, {\bf x}, {\bf w}) \J{{\bf w}}^{-2} \\ &= 
\J{r}^{-2} \left( 
D_x^2 + 2y D_x + y^2 + D_y^2 +x
\right)  K(t, {\bf x}, {\bf w}) \J{{\bf w}}^{-2} 
\\ &= \sum_{k=1}^7 \CAL{Q}_{1,k}(t,{\bf x}) K(t, {\bf x}, {\bf w}) \CAL{Q}_{2,k}(t, {\bf w}),
\end{align*}
where for $k = 1,...,7$ we have
$$
Q_{1,k} (t, {\bf x}) = q_{1,k}(t) m_{1,k}({\bf x}),
Q_{2, k} (t, {\bf w}) = q_{2, k}(t) m_{2, k}({\bf w})$$ 
and 
$$|q_{1,k}(t)q_{2, k}(t)| \leq C \Bigl(1+ |\cot t| (1+ t) + |\cot t|^2(1+ t + t^2) \Bigr ),$$
$$\|m_{1, k}({\bf x})\|_{L^{\infty}_{\bf x}}\leq C, \:\|m_{2, k}({\bf w})\|_{L^{\infty}_{\bf w}}\leq C.$$
Hence we have a smoothing effect 

\begin{align}
 \left\| K_3 \phi \right\|_{L^2(\R ^2)}  \nn &= 
C  \left\| 
\int_{[0, \infty) \setminus \Omega} (\sin t)^{-1} e^{-it(-\lambda - i \nu)} \int_{\R ^2} \sum_{k=1}^7\CAL{Q}_{1,k} (t,{\bf x}) K(t,{\bf x}, {\bf w}) \CAL{Q}_{2,k}(t, {\bf w}) \phi ({\bf w}) d {\bf w} dt \right\| _{L^2(\R ^2)}
\\ \nn &\leq 
C \int_{[0, \infty) \setminus \Omega}  \left\| \sum_{k=1}^7 \CAL{Q}_{1,k} (t,{\bf x}) \left( e^{-it(H_0 - \lambda -i \nu)} \CAL{Q}_{2,k}(t, \cdot) \phi(\cdot) \right)  ({\bf x} )  \right\|_{L^2(\R _{\bf x}^2)}dt 
\\ \nn & 
\leq C \int_{[0, \infty) \setminus \Omega}   \sum_{k=1}^7 \left\| \CAL{Q}_{1,k}(t,{\bf x})  \right\|_{L^{\infty}_{\bf x}} \left\| \left( e^{-it(H_0 - \lambda -i \nu)} \CAL{Q}_{2,k}(t, \cdot) \phi (\cdot) \right)  ({\bf x} )  \right\|_{L^2(\R _{\bf x}^2)} dt  
\\ \nn & 
 \leq 
C \int_{[0, \infty) \setminus \Omega}\sum_{k=1}^7   \| \CAL{Q}_{1,k}(t, \cdot) \|_{L^{\infty}} \|  \CAL{Q}_{2,k} (t, \cdot) \|_{L^{\infty}} \|\phi \|_{L^2(\R ^2)} e^{- \nu t} dt  
\\ \nn & \leq 
C \sum_{n \in {\N}} \int_{n\pi + \lambda^{- \theta}}^{(n+1) \pi - \lambda^{- \theta}}  e^{-\nu t} \left( 
1 + \lambda^{\theta} (1 + t) +  (\sin t)^{-2} (1+ t + t^2)
\right)  \| \phi \|_{L^2(\R ^2)} dt 
\\ &\leq C \nu^{-1}   \lambda^ {\theta} (1 + \nu^{-1} + \nu^{-2}) \| \phi \|_{L^2(\R^2)}. \label{aX3}
\end{align}

Here we used an integration by parts for the term involving $(\sin t)^{-2}$ combined with the fact that for $t \in (n\pi + \lambda^{-\theta} , (n+1) \pi - \lambda^{-\theta} ) $ and $\lambda^{-\theta} \leq \frac{\pi}{6}$ one has a lower bound
\begin{align*}
| \sin t | = |\sin (t-n \pi)| \geq |\sin (\lambda^{- \theta})| \geq \frac{\lambda^{- \theta}}{2}.
\end{align*} 
Taking together \eqref{aX1}, \eqref{aX2} and \eqref{aX3}, we get 
\begin{align*}
\left\| 
 M_{2} (\lambda, \nu ) 
\right\|_{L^2 \to L^2}  \leq C(\theta) \nu^{-1}  \left( 
\lambda^{- \theta} + ( 1 +\nu)\lambda^{ -1} +( 1 +  \nu^{-2})  \lambda^{ \theta -1}
\right) . 
\end{align*} 
Clearly  
\begin{align*}
\left\| 
 M_{0} (\lambda, \nu )
\right\|_{L^2 \to L^2}  \leq  \nu^{-1}
\end{align*} 
 and by Lemma \ref{LL1} with $A=B= \J{r}^{-1}$, $T= (H_0 - \lambda -i \nu)^{-1}$, $\alpha _0 = \beta _0 =2$, $\alpha_1 = \beta _1 = 0$ and $\theta ' = 1- \delta /2$, one deduces
\begin{align*}
\left\| 
 M_{\delta} (\lambda, \nu ) 
\right\|_{L^2\to L^2}  \leq (C(\theta))^{\delta/2}  \nu^{-1}\left( 
\lambda^{- \theta} + (1 + \nu)\lambda^{-1} +( 1 + \nu^{-2} ) \lambda^ {\theta -1}
\right)^{\delta /2}. 
\end{align*}
\end{proof}

\section{Absence of large embedded eigenvalues}
\renewcommand{\theequation}{\arabic{section}.\arabic{equation}}
\setcounter{equation}{0} 

In this section we study the relation
$$\sigma_{pp}(H) \cap \Bigl((-\infty, -R_1)\cup(R_2, \infty)\Bigr) = \emptyset.$$ 
and we work without any assumption on the support of $V$. 
The absence of large eigenvalues has been established by Dimassi-Petkov \cite{DP}. However, the fact that $R_1, R_2 > 0$ do not  depend on the support of $V$ has not been proven in \cite{DP}.  Here we establish this result and, moreover, we obtain bounds for $R_1, R_2.$

In Appendix B we prove the following
\begin{Prop}\label{eq:4.1}
Assume that we have 
\begin{align}
\|\la x \ra^{1/2} \la y \ra V\|_{L^{\infty}} \leq A_2,  
\la x \ra^{1/2} V \to 0,\:\la x \ra^{1/2} V_x \to 0, \:\la y \ra V_x \to 0 \:\: {\rm as}\: (x^2 + y^2) \to \infty. \label{eq:4.1} 
\end{align} 
Let $\psi$ be an eigenfunction of $H$ with eigenvalues $\lambda$. Then $D_x\psi\in L^2(\R^2)$.
\end{Prop}

\begin{Prop}\label{P1}
Assume that $V$ satisfies the conditions $(\ref{eq:4.1})$ and
$$\sup_{(x, y) \in \R^2} |\J{x} \J{y} V_x(x,y) | \leq A_1.$$
 Then there exist  constants $R_1 > 0, R_2 >0$ independent of $\eta_0$ such that 
 $$\sigma_{pp}(H) \cap \Bigl((-\infty, -R_1)\cup(R_2, \infty)\Bigr) = \emptyset.$$
 Moreover, we have 
 \begin{equation}\label{eq:4.2}
 R_2 \leq (C_a A_1)^8,\: R_1 \leq C_a + A_0,
 \end{equation} 
 where $C_a > 0$ is a constant depending on the choice of the function $\varphi(x)$ in Section $2$ and $a > 0.$
\end{Prop}
Notice that in the case when $V$ satisfies Assumption 1.1 the conditions of Proposition 4.2 are fulfilled.
\begin{proof} Let $\psi$ and $\lambda$ be an eigenfunction and an eigenvalue of $H$, respectively. Let $\|\psi\| = 1$ and let $|\lambda| \geq 1$. The operator $D_x$ is a conjugated operator for $H$ in the sense of \cite{M} and $D_x$ satisfies the conditions (a)-(e) in \cite{M} (see for more details Section 3 in \cite{DP2}). In particular, the condition (c) in \cite{M} means that for $\Psi \in \D{H} \cap \D{D_x} $ the symmetric form
$$(\Psi, i [H, D_x]  \Psi ) =  i (H\Phi, D_x \Psi) - i (D_x \Psi, H \Psi)$$
is bounded from below and closable and we can define the self-adjoint operator $[H, D_x]^o$ associated to its closure (\cite{M}). According to Proposition 4.1, we have $\psi \in \D{H} \cap \D{D_x}$. Thus $i[H , D_x]^o \psi$ is well defined and  $0 = (\psi, i[H, D_x]^o \psi) = ((1 + V_x)\psi, \psi).$ Consequently,

\begin{align} \nn
0 =\Bigl |\left( 
i[H, D_x]^o \psi , \psi
\right)\Bigr | & \geq 1 - \Bigl|\left( 
(\partial _x V) \psi , \psi
\right) \Bigr| \\
& \geq  1 -  \left\|
 \J{x} \J{y} \pa_x V \right\|_{L^{\infty}} \left\| \J{x}^{-1} \J{y}^{-1}  \psi \right\|.  \label{eq:4.3}
\end{align}

\def\phi{\varphi}
\def\epsilon{\varepsilon}
\def\ii{{\bf i}}
\def\Rc{{\mathcal R}}
\def\tE{\tilde{E}}
\def\ep{\epsilon}
\def\la{\langle}
\def\ra{\rangle}
\def\co{{\mathcal O}}
\def\pa{\partial}

 Let $\varphi(x)$ be the function introduced in Section 2. Obviously, with a constant $c_a > 0$ one has
\begin{equation} \label{eq:4.4}
|\varphi''(x) (\varphi' (x))^{-1}  | \leq c_a, \: \forall x \in \R.
\end{equation} 
Recall that from Lemma 2.6 we have
\begin{equation}\label{eq:4.5}
\|\svx  (D_x + y) \psi\| \leq C_1 |\lambda|^{1/4},\: \|\svd D_y \psi\|\leq C_1 |\lambda|^{3/8}.
\end{equation}

We need the following
\begin{Lem} We have the equality
\begin{align}\label{eq:4.6} 
\Gamma(\psi) &:=\|\svd (D_x + y) \psi\|^2 + \|\svd D_y \psi\|^2 \nn \\
& \quad - \Im  \Bigl( \varphi''(x) (\varphi'(x))^{-1} \svd (D_x + y) \psi , \svd \psi\Bigr)\nn \\
& \qquad - \Im \Bigl(( \varphi''(y) (\varphi'(y))^{-1} \svd D_y \psi , \svd \psi\Bigr)\nn \\
& \quad \qquad +  \Bigl(\svd (x + V) \psi, \svd \psi\Bigr) 
\\ &= \lambda (\psi, \varphi'(x) \varphi'(y)\psi). \nonumber
\end{align} 
\end{Lem} 
Notice that by (\ref{eq:4.5}) all scalar products in (\ref{eq:4.6}) are well defined.
\begin{proof}
Choose a sequence of functions $f_n \in C_0^{\infty}(\R^2)$ such that $f_n \to \psi, Hf_n \to H\psi$ in $L^2.$
Clearly,
 $$((H- \lambda) f_n, \varphi'(x) \varphi'(y) f_n) \rightarrow ((H- \lambda)\psi, \varphi'(x) \varphi(y) \psi) = 0.$$
 By integration by parts, we will show that
 \begin{equation} \label{eq:4.7}
 \Gamma(f_n) = (Hf_n, \varphi'(x) \varphi'(y) f_n)
 \end{equation}
 which yields 
 $$((H - \lambda)f_n, \varphi'(x) 
 \varphi'(y)  f_n) = \Gamma(f_n) - \lambda(f_n, \varphi'(x) \varphi'(y)  f_n).$$

To do this, we  transform the term
$$((D_x + y)^2 f_n + D_y^2 f_n, \varphi'(x) \varphi'(y) f_n).$$
First consider
\begin{align*}
(D_y^2 f_n, \varphi'(x) \varphi'(y) f_n) = (\varphi'(y)  D_y \svx D_y f_n, \sqrt{\varphi'(x)} f_n)\nonumber \\
 = (\svd D_y f_n, \svd D_y f_n)\\
 + i( \varphi''(y) (\varphi'(y))^{-1} \svd D_yf_n, \svd f_n).
 \end{align*}
Second, by the same argument we get
\begin{align*}
((D_x + y)^2 f_n, \varphi'(x) \varphi'(y) f_n) =\Bigl (\varphi'(x) (D_x + y) \svy (D_x + y) f_n, \svy f_n\Bigr) \nn\\
=( \svd (D_x + y) f_n, \svd (D_x + y) f_n)\\
+ i( \varphi''(x) (\varphi'(x))^{-1} \svd (D_x + y) f_n , \svd f_n).
\end{align*} 
Thus we obtain (\ref{eq:4.7}).
We take the limit $n \to \infty$ and deduce $\Gamma(f_n ) \rightarrow \Gamma(\psi).$ Indeed, by Lemma 2.2 we have in $L^2$ the convergence 
$$\svx (D_x + y) f_n \to \svx (D_x + y) \psi,\: \svx D_y f_n \to \svx D_y \psi$$ 
and the function $\varphi'(x) x$ is bounded for all $x \in \R.$
\end{proof}
Applying (\ref{eq:4.4}) and (\ref{eq:4.5}), one has
$$\Bigl| \Bigl( \varphi''(y) (\varphi'(y))^{-1} \svd D_y \psi, \svd \psi\Bigr)\Bigr| $$\
$$\leq C \|\sqrt{ \varphi(x)\varphi'(y)} D_y\psi\|
\leq C|\lambda|^{3/8},$$
$$\Bigl| ( \varphi''(x) (\varphi'(x))^{-1} \svd (D_x +y)  \psi, \svd \psi)\Bigr| $$
$$\leq C\|\sqrt{ \varphi(x)\varphi'(y)} (D_x + y) \psi\|
 \leq C |\lambda|^{1/4}.$$
Consequently, from (\ref{eq:4.6})  one deduces
$$\|\svd \psi\|^2 \leq C (|\lambda|^{3/4} + 1) |\lambda|^{-1} \leq C|\lambda|^{-1/4}, \: |\lambda| \geq 1,$$
hence
$$\| \la x \ra^{-1} \la y\ra^{-1} \psi\| \leq C_0 \|\svd \psi\| \leq C_1 |\lambda|^{-1/8}.$$
Going back to (\ref{eq:4.3}), we deduce that for $ |\lambda| \geq (2 C_1 A_1)^8$ we have no eigenvalues of $H$.\\

   For $\lambda \leq 0$ we have better result. For simplicity of notations denote
   $$\|\svd (D_x + y) \psi\| = B_1, \: \|\svd D_y \psi\| = B_2, \|\svd \psi\| = D.$$
   Since $-\lambda \|\svd \psi\| \geq 0,$ the equality  (\ref{eq:4.6}) implies
$$ B_1^2 + B_2^2 - C_2 B_1 D - C_3 B_2 D  + \Bigl(\svd (x + V) \psi, \svd \psi\Bigr)  \leq  - |\lambda| D^2$$
with constants $C_2 > 0, C_3 > 0$ independent of $\lambda$. Therefore,
$$\Bigl (B_1 - \frac{C_2}{2} D\Bigr)^2  + \Bigl(B_2 - \frac{C_3}{2} D\Bigr)^2 - \Bigl(\frac{C_2^2}{4} + \frac{C_3}{4}\Bigr)  D^2 \leq (C_4 + A_0 - |\lambda|)  D^2$$
with a constant $C_4 > 0$ depending on $\varphi(x) $ and independent of $\lambda.$
Consequently,  one deduces
$$ |\lambda| D^2 \leq \Bigl(\frac{C_2^2}{4} + \frac{C_3}{4}\Bigr)  D^2 + (C_4 + A_0) D^2 = (C_5 + A_0) D^2$$
If $|\lambda| > C_5 + A_0$, we have
$$ \|\varphi'(x) \varphi'(y) \psi\| = 0,$$
hence $\psi = 0.$
\end{proof}

\section{Mourre type estimate for the operator $H_0$}
\renewcommand{\theequation}{\arabic{section}.\arabic{equation}}
\setcounter{equation}{0}

In this section we fix $R \geq \max\{ R_1, R_2\}$, where  $R_k, \: k =1, 2,$ are given by Proposition 4.2.  The following result follows from \cite{M}.
 \begin{Prop}\label{PPP1}
 There exists a constant $C_R > 0$ such that 
\begin{align} \label{eq:5.1}
& \sup_{\lambda \in [-R,R] , \nu >0} \left\| 
|D_x + \beta + i|^{-1} (H_0 - \lambda  \mp i \nu)^{-1} |D_x + \beta + i|^{-1}
\right\|_{L^2 \to L^2}  \leq C_R.
\end{align} 
\end{Prop}  
 We have
$$i[D_x +  \beta, H_0] = 1.$$
As it was mentioned in the previous section,  the conjugate operator $D_x + \beta$ satisfies the conditions (a)-(e) in \cite{M} and  the principal theorem in \cite{M} implies the  estimate (\ref{eq:5.1}).

\begin{Prop}\label{P3}
Let $0 < \gamma < 1/2, \: s = 1/2 + \gamma/2,\: \beta \in \R$ and $\lambda \in [-R,R]$. Then we have the estimate 
\begin{align} \label{eq:5.2}
\sup_{\lambda \in [-R,R], \nu >0} \left\| |y- \beta|^{-\gamma} F_{\eta _0}(y - \beta)\J{x}^{-s}
 (H_0-\lambda \mp i \nu )^{-s}\J{x}^{-s}F_{\eta _0}(y - \beta) |y- \beta|^{-\gamma}
\right\|_{L^2 \to L^2}  \leq C_{R, \gamma} 
\end{align}
with a constant  $C_{R, \gamma}>0$  independent of $\eta_0$ and $\beta$.
\end{Prop} 
\begin{proof}For simplicity we treat the case $\beta = 0$. Define 
\begin{align*}
\CAL{J}_{H_0} :=  |D_x + i|^{-1}
 (H_0-\lambda \mp i \nu )^{-1} |D_x + i|^{-1}.
\end{align*} We write 
\begin{align*}
& |y|^{-\gamma} F_{\eta _0}(y ) \J{x}^{-s}(H_0-\lambda \mp i \nu )^{-1}  \J{x}^{-s} F_{ \eta _0} (y ) |y|^{- \gamma} 
\\ &= 
 |y|^{-\gamma} F_{\eta _0}(y ) \J{x}^{-s}  F_{2R}(H_0)  (H_0-\lambda \mp i \nu )^{-1} F_{2R}(H_0) \J{x}^{-s}  F_{\eta _0} (y ) |y|^{- \gamma} \\ & \qquad + |y|^{-\gamma} F_{ \eta _0}(y ) \J{x}^{-s} (1-F_{2R}^2(H_0))  (H_0-\lambda \mp i \nu )^{-1} \J{x}^{-s}  F_{\eta _0} (y ) |y|^{- \gamma} 
\\ &=
I_1 \CAL{J}_{H_0} I_1^{\ast}  +  I_2
\end{align*}
with 
\begin{align*}
I_1:= |y|^{-\gamma} F_{\eta _0}(y ) \J{x}^{-s} F_{2R}(H_0) |D_x + i|
\end{align*}
and 
\begin{align*}
I_2:= |y|^{-\gamma} F_{\eta _0}(y )\J{x}^{-s} (1-F_{2R}^2(H_0)) (H_0-\lambda \mp i \nu )^{-1} \J{x}^{-s} F_{\eta _0} (y ) |y|^{- \gamma}  .
\end{align*}
First, we show that $\left\| I_1 \right\|_{L^2 \to L^2}  \leq C_{1, R,\gamma}.$

To do this, one considers the product
\begin{align*}
I_1 =  I_{1,1} I_{1,2} 
\end{align*} 
with 
\begin{align*} 
I_{1,1} &:= |y|^{- \gamma} F_{\eta _0} (y) \J{x}^{-s} F_{2R} (H_0)  (D_x + i) , \\ 
I_{1,2} &:= (D_x +i)^{-1}|D_x  + i|. 
\end{align*}
Clearly, $I_{1,2} $ is a bounded operator. \\

  Next we write
\begin{align*} 
I_{1,1} &=  |y|^{- \gamma} F_{\eta _0} (y) \J{x}^{-s} (D_x + i) F_{2R} (H_0) \\ & \quad + |y|^{- \gamma} F_{\eta _0} (y) \J{x}^{-s}\J{H_0}^{-\gamma} \J{H_0}^{\gamma}  [F_{2R}(H_0), D_x] = J_1 + J_2.
\end{align*}
The term $J_1$ can be estimated by
\begin{align*}
& \left\| 
|y|^{-\gamma} F_{1} (y) \J{x}^{-s} (D_x + y + i) F_{2R}(H_0)
\right\| + \left\| 
 |y|^{1-\gamma} F_{1} (y) \J{x}^{-s} F_{2R}(H_0)
\right\| \\ & 
\leq 
C \left\| 
\J{D_y}^{\gamma} F_{1} (y) \J{x}^{-s} (D_x + y + i) F_{2R}(H_0)
\right\|  + C_1 \\ & \leq 
 C \left\| 
\J{D_y}^{\gamma} F_{1} (y) \J{x}^{-s} (D_x + y + i) (H_0 + i)^{-2}
\right\|   + C_1.
\end{align*} 
To handle the operator on the right hand side, write
\begin{align} 
& \nn  \J{D_y}^{\gamma}  F_{1} (y) \J{x}^{-s} (D_x + y + i) (H_0 +i)^{-2} 
\\ & \label{eq:5.3} = 
  \J{D_y}^{\gamma}  F_{1}(y)  \J{x}^{-\gamma/2} (H_0+i)^{-1}  \J{x}^{-1/2}  (D_x + y + i) (H_0 +i)^{-1} 
  \\ & \nn \quad + 
   \J{D_y}^{\gamma}   F_{1}(y)  \J{x}^{-\gamma/2} [    \J{x}^{-1/2}  (D_x + y + i) , (H_0 + i)^{-1}] (H_0 +i)^{-1}. 
\end{align}
According to  \eqref{eq:2.2}, \eqref{eq:2.3} and Corollary 2.6,  the first term in right hans side of  \eqref{eq:5.3} is bounded. For the second term one has
\begin{align*}
& \J{D_y}^{\gamma}   F_{1}(y)  \J{x}^{-\gamma/2} [    \J{x}^{-1/2}  (D_x + y + i) , (H_0 + i)^{-1}] (H_0 +i)^{-1} \\ &= 
\J{D_y}^{\gamma}   F_{1}(y) \J{x}^{-\gamma/2} (H_0 + i)^{-1} [ H_0 ,   \J{x}^{-1/2}  (D_x + y + i) ] (H_0 +i)^{-2} .
\end{align*}
Clearly, 
\begin{align*}
& [ H_0 ,   \J{x}^{-1/2}  (D_x + y + i) ] (H_0 +i)^{-2}  \\ & =[ (D_x + y)^2 + D_y^2 + x ,   \J{x}^{-1/2}  (D_x + y) + i \J{x}^{-1/2} ] (H_0 +i)^{-2} \\ &= 
i  x \J{x}^{-5/2} (D_x + y)^2 (H_0 + i)^{-2} /2 + \CAL{B}_0, 
\end{align*}
where $\CAL{B}_0$ is a bounded operator. It remains to show that the operator
$$\CAL{B}_1= \la x \ra^{-1} (D_x + y)^2 (H_0 - i)^{-2}$$
 is bounded. Set $Q = (D_x + y)^2 + D_y^2$. 
 Then
 $$(D_x + y)^2(H_0 - i)^{-1} = (D_x + y)^2 (Q - i)^{-1} + (D_x + y)^2 (Q- i)^{-1} x (H_0 - i)^{-1}.$$
The pseudoffiferential operator $(D_x + y)^2(Q- i)^{-1}$  has symbol  in $S^0(\R^4_{(x,y,\xi,\eta)})$, hence it is bounded (see  \cite{DP2}). Consequently, the operator 
$$\J{x}^{-1} (D_x + y)^2 (Q- i)^{-1} x$$
is also bounded since by composition of pseudodifferential operators its principal symbol is  in $S^0(\R^4_{(x, y,\xi, \eta)}).$ This implies that $\CAL{B}_1$ is bounded.

To prove the boundedness of $J_2$,
let $\tilde{g}(z) \in C_0^{\infty} (\C)$ be an almost analytic continuation of $g(s)= F_{2R_0} (s)$ such that
$$ \bar{\partial}_z \tilde{g}(z) = {\mathcal O} (|{\rm Im}\: z|^N), \:\forall N \in \N.$$
Consider the representation
$$ F_{2R}(H_0) = \frac{1}{\pi}\int \bar{\partial}_z \tilde{g} (z)(H_0 - z)^{-1} L(dz),$$ 
where $L(dz)$ is the Lebesgue measure on $\C$. Therefore
\begin{align*}
i[F_{2R}(H_0), D_x] &= \frac{i}{\pi}\int \bar{\partial}_z \tilde{g} (z)[(H_0 - z)^{-1}, D_x]  L(dz)\\
&= \frac{i}{\pi}\int \bar{\partial}_z \tilde{g} (z)(H_0 - z)^{-1} [H_0, D_x] (H_0 - z)^{-1} L(dz)\\
 &= -\frac{1}{\pi}\int \bar{\partial}_z \tilde{g} (z)(H_0 - z)^{-2} L(dz). 
\end{align*}
On the other hand, the operator 
$$ |y|^{- \gamma} F_{\eta _0} (y) \J{x}^{-s}\la H_0\ra^{-\gamma} =\J{x}^{-1/2}  |y|^{- \gamma} F_{\eta _0} (y) \J{x}^{-\gamma/2}\la H_0\ra^{-\gamma}$$
is bounded applying Proposition 2.9 with $H$ replaced by $H_0$, while
$$ \frac{1}{\pi}\int \bar{\partial}_z \tilde{g} (z) \la H_0 \ra^{\gamma} (H_0 - z)^{-2} L(dz)$$
is trivially bounded.
 Combining the above estimates, one concludes that
\begin{align*}
\left\| 
I_{1} \right\|_{L^2 \to L^2}  \leq C_{1, R, \gamma}.
\end{align*}
 Concerning $I_2$, notice that for $|\lambda| \leq R$ by the spectral Theorem  the operator
\begin{align*}
 \J{H_0}^{\gamma} (1-(F_{2R} (H_0)) ^2) (H_0- \lambda \mp i \nu ) ^{-1} \J{H_0}^{\gamma} 
\end{align*}
is bounded.  Next one obtains the estimate
\begin{align*}
\| I_2 \|_{L^2 \to L^2} \leq C \Bigl \| |y|^{- \gamma}  F_{\eta _0} (y)  \J{x}^{-\gamma/2} \J{H_0}^{-\gamma} \Bigr\|_{L^2 \to L^2}^2 \leq C_{2, R, \gamma} 
\end{align*}
by applying once more Proposition 2.9. The case $\beta \neq 0$ can be treated by a similar argument.
\end{proof}

In the next section we need a modification of Proposition 5.2 when we have a product with a right factor $|D_x + i|^{-1}.$
\begin{Prop}\label{P4}
Let $0 < \gamma < 1/2, \: s = 1/2 + \gamma/2,\: \beta \in \R$. Then we have
\begin{align} \label{eq:5.4}
\sup_{\lambda \in [-R,R], \nu >0} \left\| |y- \beta|^{-\gamma} F_{\eta _0}(y - \beta) \J{x}^{-s}
 (H_0-\lambda \mp i \nu )^{-1} |D_x + i|^{-1}
\right\|_{ L^2 \to L^2}  \leq B_{R, \gamma} 
\end{align}
with constant $B_{R, \gamma}>0$  independent of $\eta_0$ and $\beta$.
\end{Prop} 

\begin{proof} We use the notations of the proof of Proposition \ref{P3}. For $\beta = 0$ one has
$$ |y|^{-\gamma} F_{\eta _0}(y)  \J{x}^{-s}
 (H_0-\lambda \mp i \nu )^{-1} |D_x + i|^{-1}= I_1 {\mathcal J}_{ H_0}    + J_1,$$
where $I_1$ and ${\mathcal J}_ {H_0}$ are the same as in the proof of Proposition 5.2  and 
$$J_1  = |y|^{-\gamma} F_{\eta _0}(y ) \J{x}^{-s} (1-F_{2R}(H_0)) (H_0-\lambda \mp i \nu )^{-1} |D_x + i|^{-1} .$$
Notice that the operator $J_1$  can be bounded by $C_{3, R, \gamma}$ by a calculation similar to that used for $I_2$ in the proof of Proposition 5.2 and we leave the details to the reader. The case $\beta \neq 0$ is treated by a similar argument.
 \end{proof}

\section{ Absence of embedded eigenvalues for potentials with small support}
\renewcommand{\theequation}{\arabic{section}.\arabic{equation}}
\setcounter{equation}{0}

In this section  we prove Theorem 1.3.

\begin{proof} 
Concerning $H_0$ and $0 < \gamma < 1 /2, \: s = 1/2 + \gamma/2 ,$  we have the estimates (\ref{eq:5.2}) with $\beta = 0$ and (\ref{eq:5.4}).
 For the operator $H = H_0 + V$ with ${\rm supp} \: V \subset \{(x,y):\: |y| \leq\eta_0\}$ write
 \begin{equation} \label{eq:6.1}
 (H - \lambda - i \nu)^{-1} = (H_0 -\lambda - i \nu)^{-1} \Bigl[ 1 - V (H - \lambda - i \nu)^{-1}\Bigr]
 \end{equation}
 which yields
 \begin{align*} & |y|^{-\gamma}F_{\eta_0}(y)  \J{x}^{-s} (H - \lambda - i \nu)^{-1}|D_x + i|^{-1} \\ &= |y|^{-\gamma} F_{\eta_0}(y)  \J{x}^{-s}(H_0 - \lambda - i \nu)^{-1} |D_x + i|^{-1}
 \\ & -\Bigl[|y|^{-\gamma} F_{\eta_0}(y)\J{x}^{-s}  (H_0 - \lambda - i \nu)^{-1} \J{x}^{-s}|y|^{-\gamma} F_{\eta_0}(y)\Bigr]\Bigl(\J{x}^{2s} |y|^{2\gamma} V\Bigr)  \\ & \qquad \times\Bigl[|y|^{-\gamma} F_{\eta_0}(y) \J{x}^{-s} (H - \lambda-i \nu)^{-1} |D_x + i|^{-1} \Bigr].
 \end{align*}
 Therefore,
 \begin{align*} 
& \Bigl(I + \Bigl[|y|^{-\gamma} F_{\eta_0}(y) \J{x}^{-s}(H_0 - \lambda - i \nu)^{-1} \J{x}^{-s}|y|^{-\gamma} F_{\eta_0}(y)\Bigr] \Bigl(\J{x}^{1 + \gamma} |y|^{2\gamma} V\Bigr) \Bigr)\\ 
& \qquad \times |y|^{-\gamma}F_{\eta_0}(y)  \J{x}^{-s}(H - \lambda - i \nu)^{-1} |D_x + i|^{-1}\\
& =|y|^{-\gamma} F_{\eta_0}(y)  \J{x}^{-s} (H_0 - \lambda - i \nu)^{-1} |D_x + i|^{-1}. 
\end{align*}
  Clearly,
 $$\Bigl\|\J{x}^{1 + \gamma} |y|^{2 \gamma} V\Bigr\|_{L^2(\R^2) \to L^2(\R^2)} \leq \eta_0^{2 \gamma} \| \J{x}^{1 + \gamma} V\|_{L^{\infty}(\R^2)}.$$
   Consequently, assuming  $\eta_0^{2 \gamma} C_{R, \gamma} \|\J{x}^{1 + \gamma}V\|_{L^{\infty}(\R^2)} = c_{R, \gamma, \eta_0}< 1$, we deduce that the operator in the brackets $\Bigl( ...\Bigr)$ is invertible and
 
 $$\sup_{|\lambda| \leq R, \nu > 0}\Bigl\| |y|^{-\gamma}F_{\eta_0}(y) \J{x}^{-s} (H - \lambda - i \nu)^{-1}  |D_x + i| ^{-1}\Bigr\|_{L^2 \to L^2}  \leq  \frac{ B_{R, \gamma}}{1 - c_{R, \gamma, \eta_0}}. $$
 This estimate implies that $H$ has no eigenvalues in $[-R, R].$ In fact, let $\psi$ be an eigenfunction of $H$ with eigenvalue $\lambda \in [-R, R].$ By Proposition 4.1 we know that
 $D_x \psi \in L^2(\R^2),$ hence $|D_x + i| \psi = |D_x + i| (D_x + i)^{-1} (D_x + i) \psi \in L^2(\R^2).$ Then we conclude that
 $$|y|  ^{-\gamma} F_{\eta_0}(y) \J{x}^{-s} (H - \lambda - i \nu)^{-1}  |D_x + i|^{-1} |D_x + i| \psi  = | y|^{-\gamma} F_{\eta_0}(y)\J{x}^{-s}  i \nu^{-1} \psi.$$
 If $F_{\eta_0}(y) \psi(x, y) = 0$, then $V(x, y) \psi(x, y) = 0$ and $\psi$ will be an eigenfunction of $H_0$ which is impossible. Thus $ |y|^{-\gamma} F_{\eta_0}(y) \J{x}^{-s} \psi \neq 0$ and  as $\nu \searrow 0$ the $L^2(\R^2)$ norm of the function  $|y| ^{-\gamma} F_{\eta_0}(y) \J{x}^{-s}  \nu^{-1} \psi$ is not bounded. We obtain a contradiction and the proof is complete. \end{proof}

\renewcommand{\theLem}{A.\arabic{Lem}}
 \renewcommand{\theequation}{A.\arabic{equation}}
  \setcounter{equation}{0} 
       \section*{ Appendix A } 
       \def\thesection{A\Alph{section}}

We prove in this Appendix the following\\

{\bf Lemma A.1.} 
{\it The operators 
\begin{align*}
 x (H_0-i)^{-1}, \quad (D_x + y)^k (H_0 -i)^{-1}, \quad (D_y)^k(H_0 -i)^{-1},\: k = 1,2
\end{align*}
are unbounded  from  $L^2(\R^2)$ into $L^2(\R^2)$.}  
\begin{proof} 
Set $U_1= e^{iD_xD_y}$. We have 
$$U_1^{-1}  (D_x+y) U_1 = y, \:\:U_1^{-1} x U_1 =x-D_y .$$
Combining this with the fact that $U_1$ commutes with $D_y$, we get 
$$U_1^{-1} (D_x+y)^kU_1 U_1^{-1}\left((D_x+y)^2+D_y^2+x-i\right)^{-1}U_1=y^k\left(y^2+D_y^2+x-D_y-i\right)$$
$$=y^k\left(y^2+(D_y-\frac{1}{2})^2+x-\frac{1}{4}-i\right)^{-1}.$$
Hence, applying  the unitary transformation $e^{iy/2}$, one deduces that $(D_x + y)^k (H_0 +i)^{-1}$ is unitarily equivalent to
$$L_k=y^k\left(y^2+D_y^2+x-\frac{1}{4}-i\right)^{-1}=:y^k(B-i)^{-1}.$$
Next, we prove that $L_1$ is unbounded from $L^2(\mathbb R^2)$ into $L^2(\mathbb R^2).$ Let $\varphi\in C^\infty_0(]1,2[;\mathbb R)$ be a function such that  $\int \varphi(x)^2dx=1$, and let $\psi_n(y)$ be the normalized eigenfunction  of the harmonic oscillator corresponding to $\lambda_n=2n+1$, that is
\begin{equation}\label{OS}
(D_y^2+y^2)\psi_n(y)=(2n+1)\psi_n(y), \Vert \psi_n\Vert=1.
\end{equation}
Set $\Psi_n(x,y)=\psi_n(y)\varphi(x+2n+1)$.  Clearly,
$$L_1\Psi_n(x,y)=y\left(2n+1+x-\frac{1}{4}- i \right)^{-1} \Psi_n(x,y),\,\, \Vert \Psi_n\Vert=1.$$
Therefore,
$$\Vert L_1\Psi_n\Vert^2=\int_{\mathbb R}y^{2}\psi_n^2(y)dy  \int_{\mathbb R}  \frac{\varphi^2(x+2n+1)} {(x+2n+\frac{3}{4})^2+1}dx.$$
 On the support of $\varphi(x+2n+1)$ we have $\frac{3}{4}  \leq x + 2n + \frac{3}{4} \leq 1 + \frac{3}{4}$, hence 
 $$\frac{1} {(x+2n+\frac{3}{4})^2+1}\geq \frac{16}{65}.$$
 This yields
\begin{equation}\label{OS2}
\Vert L_1\Psi_n\Vert^2 \geq \frac{16}{65}\int_{\mathbb R}y^{2}\psi_n^2(y)dy\int_\mathbb R \varphi^2(x+2n+1)dx= \frac{16}{65}\int_{\mathbb R}y^{2}\psi_n^2(y)dy.
\end{equation}
By using the Fourier transform ${\mathcal F}_{y \to \eta}$ with respect to $y$, one obtains  ${\mathcal F}(D_y^2+y^2){\mathcal F}^{-1}=D_{\eta}^2 +\eta^2$ and 
$$ \|\psi_n(y)\|=\|\hat{\psi}_n(\eta)\|.$$
Thus we deduce that $\hat{\psi}_n(\eta)$ is also a solution of \eqref{OS} and $\hat{\psi}_n(\eta) = \psi(\eta)$. Therefore
 $$\|\psi_n'(y)\|=\|\eta \psi_n(\eta)\|=\|y \psi_n(y)\|.$$
Combining this with the obvious equality
$$2n+1=\langle (D_y^2+y^2)\psi_n,\psi_n\rangle=\Vert \psi_n'\Vert^2+\Vert y\psi_n\Vert^2,$$
we deduce  that $\Vert y\psi_n\Vert^2=\frac{2n+1}{2}$. Consequently,  \eqref{OS2}  yields
$$\Vert L_1\Psi_n\Vert^2\geq \frac{16(2n+1)}{130}.$$
Letting $n\rightarrow \infty$, we conclude that $L_1$ is unbounded from $L^2$ into $L^2$. On the other hand, from
$$\Vert L_1 u\Vert^2=\vert \langle L_2u, (B-i)^{-1}u\rangle\vert\leq \Vert L_2u\Vert \Vert u\Vert,$$
we deduce that $L_2$ is also unbounded. This shows that $(D_x+y)^k(H_0-i)^{-1}$, $k=1,2$ are unbounded. Similar arguments show that $ x (H_0-i)^{-1}$ and
$(D_y)^k(H_0 -i)^{-1}
$ are unbounded.
\end{proof}

 \renewcommand{\theequation}{B.\arabic{equation}}
  \setcounter{equation}{0}  
  \section*{Appendix B} 
         \def\thesection{B\Alph{section}}

In this Appendix we establish Proposition 4.1.

\def\D {\la D_x \ra}
\def\vep{\varphi_{\ep}}
\def\tep{\tilde{\varphi}_{\ep}}
\def\ct{\chi_t(D_x)}
\def\he{h_{\ep}}
\def\He{H_{\ep}}

{\it Proof of Proposition $4.1$} 
Let $\psi $ be a normalized by $\|\psi\| = 1.$ Suppose that
\def\vep{\varphi_{\ep}}

\begin{equation}\label{eq:B1}
D_x\psi\not\in L^2(\R^2)
\end{equation}
and for $\epsilon>0$, introduce the function  $f_\epsilon(x)=\ln\left(\frac{\langle x\rangle}{1+\epsilon \langle x\rangle}\right)$. The operators $F_\epsilon=e^{f_\epsilon(D_x)}=\frac{\langle D_x\rangle}{1+\epsilon \langle D_x\rangle}$ and its inverse $F_\epsilon^{-1}=e^{-f_\epsilon(D_x)}=\frac{1+\epsilon \langle D_x\rangle} {\langle D_x\rangle}$  are  bounded. Therefore, $F_\epsilon \psi\in L^2(\R^2)$.  The condition \eqref{eq:B1} implies $\lim_{\epsilon \searrow 0}\Vert F_\epsilon \psi\Vert=\infty$.  
Let ${\mathcal F}_x = {\mathcal F}_{x \to \xi} $ denotes the Fourier transform with respect to $x$. 
The dominated convergence theorem yields 
\begin{align*} 
\lim_{\epsilon \searrow 0}\iint_{\R ^2} e^{f_\epsilon(D_x)}\psi(x,y) \overline {{\mathcal F}^{-1} _x g(\xi,y)}dxdy
&=
\lim_{\epsilon \searrow 0}\iint_{\R ^2} e^{f_\epsilon(\xi)}({\mathcal F_x}\psi)(\xi,y) \overline {g(\xi,y)}d\xi dy
\\ &=\iint_{\R^2}  \langle \xi\rangle({\mathcal F_x}\psi)(\xi,y) \overline {g(\xi,y)}d\xi dy,
\end{align*}
for all $g(\xi, y)\in {\mathcal F}_x(C^\infty_0(\R^2))$.  This implies 
$$\lim_{\epsilon \searrow 0}\iint_{\R ^2} e^{f_\epsilon(D_x)}\psi(x,y) \overline {h(x, y)}\frac{dxdy}{\Vert F_\epsilon \psi \Vert}=0, \quad  \forall h\in C^\infty_0(\R^2).$$ 
Consequently,  the normalized function $\varphi_{\epsilon}:=\frac{F_\epsilon \psi}{\Vert F_\epsilon \psi\Vert}$ converges weakly to zero.

By using  $F_\epsilon^{-1} x F_\epsilon= x+i(\partial_x f_\epsilon)(D_x)$,   and taking into account that
 $F_\epsilon$ commutes with the operator $H-x-V$, we get
\begin{equation}\label{eq:B2}
H \varphi_\epsilon=F_\epsilon \left(H+i(\partial_x f_\epsilon)(D_x)-V+F_\epsilon^{-1} V F_\epsilon\right)\frac{\psi}{\Vert F_\epsilon \psi\Vert}=(\lambda+i(\partial_xf_\epsilon)(D_x)+V-F_\epsilon V F_\epsilon^{-1})\varphi_\epsilon.
\end{equation}
Notice that the operators $\frac{1}{1+\epsilon\langle D_x\rangle}$ and $\frac{\epsilon \langle D_x\rangle}{1+\epsilon \langle D_x\rangle}$ are bounded  from $L^2$ into $L^2$ uniformly with respect to $\epsilon\in [0,1]$. Since  $V, \partial_xV\in L^\infty(\mathbb R^2)$,
 the operator $\langle D_x\rangle V\langle D_x\rangle^{-1}$ is bounded. Hence,
\begin{equation}\label{eq:B3}
F_\epsilon V F_\epsilon^{-1}=\frac{\langle D_x\rangle}{1+\epsilon\langle D_x\rangle} V\frac{1+\epsilon\langle D_x\rangle}{\langle D_x\rangle}
=\frac{1}{1+\epsilon\langle D_x\rangle}\left(\langle D_x\rangle  V \langle D_x\rangle^{-1}\right) +\frac{\epsilon\langle D_x\rangle}{1+\epsilon\langle D_x\rangle} V,
\end{equation}
is uniformly bounded for $\epsilon \in [0,1]$.

From now on we denote
$$K_\epsilon:=i(\partial_xf_\epsilon)(D_x)+V-F_\epsilon V F_\epsilon^{-1}.$$
 Let $G(x,y)$ be a continuous function going to zero as $(x^2+ y^2) \to \infty$. It is well known that $$ \phi \langle D_x\rangle^{-s} (H+i)^{-1},$$ is a compact operator for every $\phi\in C^\infty_0(\mathbb R^2)$ and all $s\geq 0$ (see for instance, \cite{DP2}). Thus,
by an approximation argument, $G\langle D_x \rangle^{-s} (H+i)^{-1}$ is also  compact. We claim that
\begin{equation}\label{eq:B4}
G\langle D_x\rangle^{-s}\varphi_\epsilon  \text { converges strongly to zero as } \epsilon \searrow 0.
\end{equation}
To prove this, we use \eqref{eq:B2}. Write
$$G\langle D_x\rangle^{-s}\varphi_\epsilon=G\langle D_x\rangle^{-s}(H+i)^{-1} (H+i)\varphi_\epsilon=G \langle D_x\rangle^{-s}
(H+i)^{-1}\left(\lambda+i + K_\epsilon\right)\varphi_\epsilon.$$
Since  $\left(\lambda+i+ K_\epsilon\right)$ is bounded uniformly for $\epsilon\in [0,1]$,  and $\varphi_\epsilon$ converges weakly to zero, it follows from the compactness of $G\langle D_x\rangle^{-s}(H+i)^{-1}$ that the right hand side of the above equality converges strongly to zero. 

For $t >1,$ let  $\chi_t(x)$ be  an odd smooth function satisfying
\begin{equation}\label{eq:B5}
\chi_t(x)=\begin{cases}  x,\: 0\leq  x  \leq t, \\[8pt]
2t,\: x\geq 2t,
\end{cases}
\end{equation}
$\chi_t^{(k)}(x)={\mathcal O}(t^{-k+1}), k\geq 1,$ and $\chi'_t(x)\geq 0$.   
Clearly,  $i[x,-\chi_t(D_x)]=\chi'_t(D_x)$ and 
\begin{equation} \label{eq:B6}
i\lim_{t \to \infty} i([x,-\chi_t(D_x)] \vep, \vep ) = (\vep, \vep).
\end{equation} 

Next, we claim that for every fixed $\ep > 0$ we have
 \begin{equation}\label{eq:B7}
 \lim_{t \to \infty}  i( [V, - \chi_t(D_x) ] \vep, \vep)=-2\lim_{t \to \infty}   \Im (\vep, \chi_t(D_x) V\vep)= (V_x \vep, \vep).
 \end{equation}
 First, it follows from \eqref{eq:B2} that  $h_{\epsilon}:=(H_0-i)\vep$ is uniformly bounded in $L^2$ with respect to $\epsilon \in [0,1]$. On the other hand, Lemma 2.3 and the conditions (\ref{eq:4.1}) show that
 $$D_xV (H_0-i)^{-1}h_\epsilon = V_x (H_0 - i)^{-1} \he + V(D_x + y)(H_0 - i)^{-1}\he - yV (H_0 - i)^{-1}\he \in L^2.$$
  Combining this with the fact that $\vert \chi_t(\xi)-\xi\vert \leq C\vert \xi\vert$ (uniformly for $t\geq 1$), we deduce 
 $$\vert (\chi_t(\xi)-\xi) H_\epsilon(\xi)\vert\leq C \vert \xi\vert\,\,\,  \vert  H_\epsilon(\xi)\vert \in L^2,  \,\,\, \text { where }\,\,\, H_\epsilon(\xi)={\mathcal F}_{x \to \xi}\Big(V (H_0-i)^{-1}h_\epsilon \Big)(\xi).$$
Hence, the dominated convergence theorem yields 
 $$\lim_{t\rightarrow +\infty} \chi_t(D_x)V(H_0-i)^{-1}h_\epsilon =V_x (H_0-i)^{-1}h_\epsilon = V_x \vep,\,\, \,\,  \text { in } \,\, L^2,$$
  and the proof of the claim is complete. Taking together \eqref{eq:B6}, \eqref{eq:B7} and the equality  $[H,\chi_t(D_x)]=[x+V,\chi_t(D_x)]$, we  obtain
  $$\lim_{t\rightarrow +\infty}i\left(\left[H,-\chi_t(D_x)\right]  \varphi_\epsilon,  \varphi_\epsilon\right)=((1+V_x)\vep,\vep).$$

Now applying \eqref{eq:B4}  with  $G=\partial_xV$ and  $s=0$, we deduce that
\begin{equation}\label{eq:B8}
 \left((1+\partial_x V\right)\varphi_\epsilon,\varphi_\epsilon) \geq \frac{1}{2} 
\end{equation}
for $\epsilon$ small enough. To complete the proof, we will show that the left hand side of  \eqref{eq:B8} is less than $\frac{1}{4}$ for
$\epsilon $ small enough. This leads to a contradiction.

Equation \eqref{eq:B2} implies
\begin{align}\label{eq:B9}
i\left(\left[H,-\chi_t(D_x)\right]  \varphi_\epsilon,  \varphi_\epsilon\right) &= i\left( \chi_t(D_x)H \varphi_\epsilon,  \varphi_\epsilon\right)-i
\left(\chi_t(D_x)\varphi_\epsilon,  H\varphi_\epsilon\right)
 \\ & \nn = i\left( \chi_t(D_x)(\lambda+K_\epsilon)\varphi_\epsilon,  \varphi_\epsilon\right)-i
\left(\chi_t(D_x)\varphi_\epsilon,  (\lambda+K_\epsilon)\varphi_\epsilon\right)
\\ & \nn =-2 \mathrm{Im}\left( \chi_t(D_x)K_\epsilon\varphi_\epsilon, \varphi_\epsilon\right).
\end{align}

On the other hand, the inequality
$$
\chi_t(x)\left(\partial_x f_\epsilon\right)(x)=\frac{ x \chi_t(x)}{\langle x\rangle^2\left(1+\epsilon \langle x\rangle\right)}\geq 0,
$$
yields
\begin{equation*}
\chi_t(D_x)\left(\partial_x f_\epsilon\right)(D_x)\geq 0,
\end{equation*} 
in the sense of self-adjoint operators. Consequently,
\begin{align}\label{eq:B10}
-2 {\text Im}\left( \chi_t(D_x)K_\epsilon\varphi_\epsilon, \varphi_\epsilon\right)\nonumber \\
=-2 {\text Im}\left(i \chi_t(D_x) (\partial_xf_\epsilon)(D_x)+V-F_\epsilon V F_\epsilon^{-1} \varphi_\epsilon, \varphi_\epsilon\right)\nonumber \\
\leq 2\:{\text Im} \left(  \chi_t(D_x)\left(F_\epsilon V F_\epsilon^{-1}-V\right)\varphi_\epsilon , \varphi_\epsilon\right).
\end{align}
From \eqref{eq:B3}, we have
$$
F_\epsilon V F_\epsilon^{-1}-V
=\frac{1}{1+\epsilon\langle D_x\rangle}\left(\langle D_x\rangle  V \langle D_x\rangle^{-1}\right) -\frac{1}{1+\epsilon\langle D_x\rangle} V 
= \frac{1}{1+\epsilon\langle D_x\rangle} \left[\langle D_x\rangle,  V\right] \langle D_x\rangle^{-1}.
$$
Therefore
\begin{align}\label{eq:B11}
\lim_{t\rightarrow \infty}\Bigl( \chi_t(D_x)\left(F_\epsilon V F_\epsilon^{-1}-V\right)\varphi_{\ep}, \varphi_{\ep} \Bigr)\nonumber &= \Bigl(D_x\left(F_\epsilon V F_\epsilon^{-1}-V\right) \varphi_{\ep}, \varphi_{\ep}\Bigr) \\ &=\Bigl(\frac{D_x}{1+\epsilon\langle D_x\rangle} \left[\langle D_x\rangle,  V\right] \langle D_x\rangle^{-1}\varphi_{\ep}, \varphi_{\ep}\Bigr)
\end{align}
 is bounded uniformly for $\epsilon \in[0,1]$. Letting $t\rightarrow \infty$, 
we deduce  from  \eqref{eq:B7}, \eqref{eq:B10} and \eqref{eq:B11}
\begin{equation} \label{eq:B12} 
\left((1+\partial_x V\right)\varphi_\epsilon,\varphi_\epsilon)\leq 2\Im \Bigl( \frac{D_x}{1+\epsilon\langle D_x\rangle} \left[\langle D_x\rangle,  V\right] \langle D_x\rangle^{-1}\varphi_\epsilon,\varphi_\epsilon\Bigr).
\end{equation}

\def\12{\frac{1}{2}}
\def\Rm#1{{\rm#1}}
\def\lx{\la x \ra}
\def\phi{\varphi}
\def\epsilon{\varepsilon}
\def\kappa{\varkappa}
\def\ep{\epsilon}

\def\Re{{\rm Re}}
\def\Im{{\rm Im}}
\def\CAL{\mathcal} 
\def\R{{\mathbb R}}
\def\N{{\mathbb N}}
\def\C{{\mathbb C}}
\def\r{\J{r}}
\def\pa{\partial}
\def\pu{\partial_t u}
\def\la{\langle}
\def\ra{\rangle}
\def\vep{\varphi_{\ep}}
\def\D{\la D_x\ra}

To complete the proof of Proposition 4.1, we apply the following\\

{\bf Lemma B.1.} 
{\it We have
$$\lim_{\ep \to 0} \Bigl(\frac{1}{1 + \ep\D} D_x[ \D, V] \D^{-1} \vep, \vep\Bigr) = 0.$$}
\begin{proof} Write
\begin{align*}
D_x[ \D, V] \D^{-1}  &=\Bigr( \D D_x V - D_x V \D \Bigr) \D^{-1} \\ & = \Bigl[ \D (V_x + V D_x) - (V_x + V D_x) \D \Bigr] \D^{-1} 
\\ &= [\D, V_x] \D^{-1} + \Bigl( \D V - V\D\Bigr)D_x \D^{-1}  
\\ &=  [\D, V_x] \D^{-1} + [\D, V] D_x \D^{-1}
\end{align*}
 and set
   $$L_1 := \frac{1}{1 + \ep \D}[\D, V_x]\D^{-1},\: L_2: = \frac{1}{1 +\ep\D}\Bigl([\D, V]D_x \D ^{-1}\Bigr).$$
Therefore, 
  $$D_x \D^{-1}  \vep = D_x \D^{-1}  (H_0- i)^{-1}(H_0 - i) \vep $$
  $$= (H_0 - i)^{-1} D_x \D^{-1} h_{\ep}  + (H_0-i)^{-1} [D_x \D^{-1}  , x] (H_0- i)^{-1} h_{\ep}.$$
  Clearly, the operator $[D_x \D^{-1}, x] = \D^{-1}( 1 - \frac{D_x^2}{\D^2})$ is bounded and this implies that
  \begin{equation} \label{eq:B13} 
  D_x \D^{-1} \vep = (H_0 - i)^{-1} \tilde{h}_{\ep}
  \end{equation}
  with $\tilde{h}_{\ep}$  bounded in $L^2$ uniformly with respect to $\ep.$ Recall that the operator
  $$\la x \ra^{-1/2}( D_x + y) (H_0 - i)^{-1} $$
  is  bounded by Lemma 2.3. By using this, one deduces that the operator
  $$(H_0 + i)^{-1} V D_x (H_0 - i)^{-1} =(H_0 + i)^{-1} V \la x \ra^{1/2} \la  x \ra^{-1/2} (D_x + y) (H_0 -i)^{-1} - (H_0 + i)^{-1} V y (H_0 - i)^{-1} $$
  is compact since $V \la x \ra^{1/2} \to 0, V y \to 0$ as $(x^2 + y^2) \to \infty$ by conditions (\ref{eq:4.1}).
 To handle the operator $\D$,  we exploit the following representation
  \begin{align*}
  V \D (H_0- i)^{-1} &= V( D_x + i) \D(D_x + i)^{-1} (H_0 - i)^{-1} \\
  & = V(D_x + i) (H_0 - i)^{-1} \D(D_x + i)^{-1} \\
  & \quad  + V(D_x + i) (H_0 - i)^{-1} [\D(D_x + i)^{-1} , x] (H_0 - i)^{-1}.
  \end{align*} 

  Obviously, the  commutator $[\D(D_x + i)^{-1}, x]$ is a bounded operator and $(D_x + i) = (D_x + y) - (y- i)$. So as above we obtain that
  $(H_0 +i)^{-1} V\D (H_0- i)^{-1} $ is compact. In the same way we show that the operators
  $$ (H_0 + i)^{-1} V_x D_x (H_0  - i)^{-1}, \: (H_0 + i)^{-1} V_x \D(H_0 - i)^{-1} $$
  are compact because $V_x \la x \ra^{1/2} \to 0, V_x y \to 0$ as $(x^2 + y^2) \to \infty$ according to conditions (\ref{eq:4.1}).

  To deal with the operator $L_1$, write 
  \begin{align*}
 & \Bigl(\frac{1} {1 + \ep \D} [\D, V_x]  \D^{-1} \vep, \vep\Bigr)
  \\ &= \Bigl(\frac{1}{1 + \ep \D} [\D, V_x] \D^{-1} \vep, (H_0 - i)^{-1} h_{\ep}\Bigr)
    \\ &= \Bigl(\frac{1}{1 + \ep\D}(H_0 + i)^{-1}  [\D, V_x] \D^{-1} \vep, h_{\ep}\Bigr) 
    \\ & \quad -\Bigl((H_0 + i)^{-1} \frac{\ep D_x}{\D(1 + \ep\D)^2} (H_0 + i)^{-1} [\D, V_x] \D^{-1} \vep, h_{\ep} \Bigr).
    \end{align*}
  We have $(H_0 + i)^{-1} [\D, V_x] \D^{-1} = (H_0 + i)^{-1} (\D V_x \D^{-1} - V_x).$ The analysis of the term with $V_x$ is easy since $(H_0 + i)^{-1}V_x$ is compact. For the other term we get 
\begin{align*}
\D^{-1}  \vep &= \D^{-1}  (H_0 -i)^{-1} h_{\ep} 
\\ & = (H_0 - i)^{-1} \D^{-1}  h_{\ep} 
  - (H_0- i)^{-1} \frac{\ep D_x} {\D(1 + \ep \D)^2} (H_0 -i)^{-1} h_{\ep}
\end{align*}  
and notice that  $(H_0 + i)^{-1} \D V_x ( H_0 - i)^{-1} $ is compact.
  
   Passing to the analysis of the operator $L_2$, we have $[\D, V] = \D V - V \D.$ For $\D VD_x \D^{-1}$ we repeat the above argument by using (\ref{eq:B13}) and the fact that 
 $(H_0 + i)^{-1} \D V(H_0 - i)^{-1} $ is compact since its adjoint $(H_0 + i)^{-1} V \D (H_0 - i)^{-1}$ is compact. On the other hand, applying (\ref{eq:B13}) once more, we have
 $$V\D D_x \D^{-1} \vep = V \D (H_0-i)^{-1} \tilde{h}_{\ep}.$$
The operator $V \D (H_0 -i)^{-1}$ has been treated above  and the proof is complete.
\end{proof}

{\bf Acknowledgment.} Thanks are due to the referees for their critical comments and useful suggestions leading to an improvement  of  the previous version of the paper.


\begin{thebibliography}{ABCD} 

\bibitem{AK} T. Adachi and M. Kawamoto, {\em Avron-Herbst type formula in crossed constant magnetic and time-dependent electric fields}, Lett. Math. Phys., \textbf{102} (2012), 65-90.

\bibitem{AD} M. Ben-Artzi and A. Devinatz, {\em Regularity and Decay of Solutions to the
Stark Evolution Equation}, J. Funct. Anal. {\textbf 154} (1998), 501-512.

\bibitem{DP1} M. Dimassi and V. Petkov, {\em Resonances for magnetic Stark hamiltonians in two dimensional case},
IMRN, \textbf{77} (2004), 4147-4179.

\bibitem{DP2} M. Dimassi and V. Petkov, {\em Spectral shift function for operators with crossed magnetic and electric fields},  Rev. Math. Physics, \textbf{22} (2010), 355-380.

\bibitem{DP} M. Dimassi, and V. Petkov, {\em Spectral problems for operators with crossed magnetic and electric fields}, J. Phys. A, \textbf{43} (2010),  474015. 

\bibitem {FK1} C. Ferrari and H. Kovarik, {\em Resonances width in crossed electric and magnetic fields}, J. Phys. A: Math. Gen. \textbf{37} (2004), 7671-7697.

\bibitem {FK2} C. Ferrari and H. Kovarik, {\em On the exponential decay of magnetic Stark resonances}, Rep. Math. Phys. \textbf{56} (2005), 197-207.

\bibitem{GS}   J. Galkowski and J. Shapiro, {\em Semiclassical resolvent bounds for weakly decaying potentials}, arXiv: math.AP: 2003.02525.

\bibitem{GM} S. Gyger, P. A. Martin, {\em Lifetimes of impurity states in crossed magnetic and electric fields}, J. Math. Phys., \textbf{40} (1999), 3275-3282. 

\bibitem{HL} E. H. Hauge, J. M. J. van Leeuwen, {\em Bound state and metastability near a scatter in crossed electromagnetic fields}, Physica A, \textbf{268} (1999), 525-552.  

\bibitem{Is} H. Isozaki, Many body Schr\"{o}dinger equation, Springer-Verlag Tokyo, (2004). (In Japanese.) 

\bibitem{Ka} M. Kawamoto,{\em Exponential decay property for eigenfunctions of Landau-Stark Hamiltonian}, Rep. Math.Phys. \textbf{77} (2016), 129-140 .

 \bibitem{M} E. Mourre, {\em Absence of singular continuous spectrum for certain self-adjoint operators}, Commun. Math. Phys.,  \textbf{78}, no. 3 (1981), 391-408. 

\bibitem{SW} E. Stein and G.  Weiss, {\em  Fractional integral on n-dimensional Euclidian space}, J. Math. Mech., \textbf{7} (1958), 03-514. 

\bibitem{Va} A. F. Vakulenko, {\em The absence of bound states for
the two-body system in the external electric field} , Zap. Nauchn. Sem. LOMI, \textbf{152} (1986), 18-20 (in Russian). English translation: J. Soviet Math. \textbf{40} (1988), 599-601.

\bibitem{V1} G. Vodev, {\em Resolvent estimates for the magnetic Schr\"odinger operator} , Anal. PDE \textbf{7}, no. 7  (2014), 1639–1648. 

\bibitem{V2} G. Vodev, {\em Semiclassical resolvent estimates for $L^{\infty}$  potentials on Riemannian manifolds}, Ann. Henri
Poincar\'e, \textbf{21} (2020), 437-459.

\bibitem{V3} G. Vodev, {\em Semiclassical resolvent estimates for H\"older potentials}, arXiv: math.AP. 2002.12853.

\bibitem{Ya} D. Yafaev, {\em Sharp constants in the Hardy-Rellich inequalities}, J. Funct. Anal. \textbf{168} (1999), 121-144.


\end{thebibliography}
\end{document}